\documentclass{scrartcl}

\usepackage{amsfonts}
\usepackage{amsmath}
 \numberwithin{equation}{section}
\usepackage{amssymb}
\usepackage{amsthm}
 \newtheorem{thm}{Theorem}[section]

 \theoremstyle{definition}
 
 \theoremstyle{remark}
 
\usepackage[english]{babel}
\usepackage[fixlanguage,noissn,noisbn]{babelbib}%
 \setbtxfallbacklanguage{english}
\usepackage{cite}%
\usepackage[T1]{fontenc}
\usepackage{hyperref}
\usepackage{lmodern}
\usepackage{mathscinet}
\usepackage{mathtools}
 \mathtoolsset{mathic=true}
\usepackage{microtype}
\usepackage{pmat}
\usepackage{url}

\usepackage{177arxiv}

\title{An Application of the Schur Complement to Truncated Matricial Power Moment Problems}
\author{Bernd Fritzsche \and Bernd Kirstein \and Conrad M\"adler}

\begin{document}

\dedication{Dedicated to the 80th birthday of M.~A.~Kaashoek}

\maketitle

\begin{abstract}
 The main goal of this paper is to reconsider a phenomenon which was treated in earlier work of the authors' on several truncated matricial moment problems.
 Using a special kind of Schur complement we obtain a more transparent insight into the nature of this phenomenon.
 In particular, a concrete general principle to describe it is obtained.
 This unifies an important aspect connected with truncated matricial moment problems.
\end{abstract}

\begin{description}
 \item[Keywords:] Truncated matricial Hamburger moment problems, truncated matricial \(\ug\)\nobreakdash-Stieltjes moment problems, Schur complement.
\end{description}

\section{Introduction}
 In this paper, we reconsider a phenomenon which was touched in our joint research with Yu.~M.~Dyukarev on truncated matricial power moment problems (see~\zita{MR2570113}, for the Hamburger case and~\zita{MR2735313} for the \(\ug\)\nobreakdash-Stieltjes case).
 
 In this introduction, we restrict our considerations to the description of the Hamburger case.
 In order to describe more concretely the central topics studied in this paper, we give some notation.
 Throughout this paper, let $p$ and $q$ be positive integers.
 Let $\N$, $\NO$, $\Z$, $\R$, and  $\C$ be the set of all positive integers, the set of all non-negative integers, the set of all integers, the set of all real numbers, and the set of all complex numbers, respectively.
 For every choice of $\rho, \kappa \in \R \cup\set{-\infty,\infp}$, let $\mn{\rho}{\kappa}\defeq \setaca{k\in\Z}{\rho\leq k \leq \kappa}$.
 We will write $\Cpq$, $\CHq$, $\Cggq$, and $\Cgq$ for the set of all complex \tpqa{matrices}, the set of all \tH{} complex \tqqa{matrices}, the set of all \tnnH{} complex \tqqa{matrices}, and the set of all \tpH{} complex \tqqa{matrices}, respectively.
 We will use $\BAR$ to denote the $\sigma$\nobreakdash-algebra of all Borel subsets of $\R$.
 For all $\Omega \in \BAR\setminus\set{\emptyset}$, let $\BA{\Omega} \defeq \BAR\cap \Omega$.
 Furthermore, we will write $\MggqO$ to designate the set of all \tnnH{} \tqqa{measures} defined on $\BAO$, \tie{}, the set of $\sigma$\nobreakdash-additive mappings $\mu\colon\BAO\to \Cggq$.
 We will use the integration theory with respect to \tnnH{} \tqqa{measures} which was worked out independently by I.~S.~Kats~\zita{MR0080280} and M.~Rosenberg~\zita{MR0163346}.
 For all $j\in \NO$, we will use $\MggquO{j}$ to denote the set of all $\sigma \in \MggqO$ such that the integral
\[
 \mpm{\sigma}{j}
 \defeq \int_\Omega x^j \sigma(\dif x)
\]
 exists.
 Obviously, if $k,\ell\in\NO$ with $k<\ell$, then it can be verified, as in the scalar case, that the inclusion $\MggquO{\ell}\subseteq\MggquO{k}$ holds true.
 Now we formulate two related versions of truncated matricial moment problems. (The Hamburger moment problem corresponds to \(\Omega=\R\).)

\begin{prob}[\mprob{\Omega}{m}{=}]
 Let \(m\in\NO\) and let \(\seqs{m}\) be a sequence of complex \tqqa{matrices}.
 Describe the set \(\MggqOsg{m}\) of all \(\sigma\in\MggquO{m}\) for which \(\mpm{\sigma}{j}=\su{j}\) for all \(j\in\mn{0}{m}\).
\end{prob}

 The just formulated moment problem is closely related to the following:

\begin{prob}[\mprob{\Omega}{m}{\leq}]
 Let \(m\in\NO\) and let \(\seqs{m}\) be a sequence of complex \tqqa{matrices}.
 Describe the set \(\MggqOskg{m}\) of all \(\sigma\in\MggquO{m}\) for which \(\su{m}-\mpm{\sigma}{m}\) is \tnnH{} and, in the case \(m\in\N\) moreover \(\mpm{\sigma}{j}=\su{j}\) for all \(j\in\mn{0}{m-1}\).
\end{prob}

\breml{H.R1032}
 If $m\in\NO $ and $\seqs{m}$ is a sequence of complex \tqqa{matrices}, then \(\MggqOsg{m}\subseteq\MggqOskg{m}\).
\erem

\breml{H.R1035}
 If $m\in\NO$ and $\seqs{m},\seqt{m}$ are two sequences of complex \tqqa{matrices} satisfying \(\su{m}-t_{m}\in\Cggq\) and \(\su{j}=t_j\) for all \(j\in\mn{0}{m-1}\), then \(\MggqOakg{t}{m}\subseteq\MggqOskg{m}\).
\erem

 In order to state a necessary and sufficient condition for the solvability of each of the above formulated moment problems in the case \(\Omega=\R\), we have to recall the notion of two types of sequences of matrices.
 If $n \in \NO$ and if $\seqs{2n}$ is a sequence of complex \tqqa{matrices}, then $\seqs{2n}$ is called \emph{\tHnnd{}} if the block Hankel matrix
\(
 \Hu{n}
 \defeq \matauuo{\su{j+k}}{j,k}{0}{n}
\)
 is \tnnH{}.
 For all $n\in\NO$, we will write  $\Hggq{2n}$ for the set of all \tHnnd{} sequences $\seqs{2n}$ of complex \tqqa{matrices}.
 Furthermore, for all $n \in \NO$, let $\Hggeq{2n}$ be the set of all sequences $\seqs{2n}$ of complex \tqqa{matrices} for which there exist complex \tqqa{matrices} $\su{2n+1}$ and $\su{2n+2}$ such that $\seqs{2(n+1)} \in\Hggq{2(n+1)}$.
 For each $n\in\NO$, the elements of the set $\Hggeq{2n}$ are called \emph{\tHnnde{}} sequences.
 Now we can characterize the situations that the mentioned problems have a solution:

\begin{thm}[\zitaa{MR1624548}{\cthm{3.2}{}}]\label{T1615}
 Let $n\in\NO$ and let $\seqs{2n}$ be a sequence of complex \tqqa{matrices}.
 Then $\MggqRskg{2n} \neq \emptyset$ if and only if $\seqs{2n} \in\Hggq{2n}$.
\end{thm}

 In addition to \rthm{T1615}, one can show that, in the case $\seqs{2n}\in\Hggq{2n}$, a distinguished molecular \tnnH{} measure belongs to $\MggqRskg{2n}$ (see,~\zitaa{MR2570113}{\cthm{4.16}{}}).

 Now we characterize the solvability of Problem~\mprob{\R}{2n}{=}.

\begin{thm}[\tcf{}~\zitaa{MR2805417}{\cthm{6.6}{}}]\label{T1708}
 Let $n\in\NO$ and let $\seqs{2n}$ be a sequence of complex \tqqa{matrices}.
 Then $\MggqRsg{2n} \neq \emptyset$ if and only if $\seqs{2n} \in\Hggeq{2n}$.
\end{thm}

 The following result is the starting point of our subsequent considerations:
 
\bthmnl{\zitaa{MR2570113}{\cthm{7.3}{}}}{T1456}
 Let \(n\in\NO\) and let \(\seqs{2n}\in\Hggq{2n}\).
 Then there exists a unique sequence \(\Reqseq{s}{2n}\in\Hggeq{2n}\) such that
\beql{T1456.A}
 \Mggoaakg{q}{\R}{\Reqseq{s}{2n}}
 =\MggqRskg{2n}.
\eeq
\ethm

 \rthm{T1456} was very essential for the considerations in~\zita{MR2570113}.
 
 Following~\zita{MR2570113}, we sketch now some essential features of the history of \rthm{T1456}.
 The existence of a sequence \(\Reqseq{s}{2n}\in\Hggeq{2n}\) satisfying \eqref{T1456.A} was already formulated by V.~A.~Bolotnikov (see~\zitaa{MR1395706}{\clem{2.12}{}} and~\zitaa{Bol}{\clem{2.12}{}}).
 However it was shown in~\zitaa{MR2570113}{\cpagestos{804}{805}} by constructing a counterexample that the proof in~\zita{MR1395706} is not correct.
 However, in the subsequent considerations of~\zitaa{MR2570113}{\csec{7}} it was then shown that the result formulated by V.~A.~Bolotnikov is correct and that the sequence \(\Reqseq{s}{2n}\) in \rthm{T1456} is unique.
 The proof of \rthm{T1456} given in~\zita{MR2570113} is constructive and does not yield a nice formula. 
 The main goal of this paper is to present a general purely matrix theoretical object which yields applied to a special case the explicit construction of the desired sequence \(\Reqseq{s}{2n}\in\Hggeq{2n}\).
 Furthermore, we will see that another application of our construction yields the answer to a similar situation connected with truncated matricial \(\ug\)\nobreakdash-Stieltjes moment problems (see \rsec{S.K}).
 
 The discovery of the above mentioned object of matrix theory was inspired by investigations of T.~Ando~\zita{And05} in the context of Schur complements and its applications to matrix inequalities.
 (It should be mentioned that Ando's view on the Schur complement is the content of \cch{5} in the book~\zita{MR2160825} which is devoted to several aspects of the Schur complement.)
 Given a \tnnH{} \tqqa{matrix} \(A\) and a linear subspace \(\cV \) of \(\Cq\), we introduce a particular \tnnH{} \tqqa{matrix} \(\SO{A}{\cV}\) which turns out to possess several extremal properties.
 Appropriate choices of \(A\) and \(\cV \) lead to the construction of the desired sequences connected to the truncated matricial moment problems under consideration.

\section{On a special kind of Schur complement}\label{S.ando}
 Against to the background of application to matrix inequalities T.~Ando~\zita{And05} presents an operator theoretic approach to Schur complements.
 
 In particular, T.~Ando generalized the notion of Schur complement of a block matrix by considering block partitions with respect to an arbitrary fixed linear subspace.
 For the case of a given \tnnH{} \tqqa{matrix} \(A\) and a fixed linear subspace \(\cV\) the construction by T.~Ando produces a \tnnH{} \tqqa{matrix} \(\SO{A}{\cV}\) having several interesting extremal properties.
 
 For our purposes it is more convenient to choose a more matrix theoretical view as used by T.~Ando~\zita{And05}.
 For this reason we use a different starting point to the main object of this section (see \rdefn{de1ab3}).
 This leads us to a self-contained approach to several results due to Ando.
 
 In the sequel, \(\Cp\) is short for \(\Coo{p}{1}\).
 Let \(\Opq\) be the zero matrix from \(\Cpq\).
 Sometimes, if the size of the zero matrix is clear from the context, we will omit the indices and write \(\NM\).
 We denote by \(\nul{A}\defeq\setaca{x\in\Cq}{Ax=\Ouu{p}{1}}\) the null space of a complex \tpqa{matrix} \(A\).
  
\breml{lm1ab3}
 If \(M \in\Cqp\) and \(\cV \) is a linear subspace of \(\Cq\), then \(\fib{M}{\cV}\defeq\setaca{x\in\Cp}{Mx\in\cV}\) is a linear subspace of \(\Cp\).
 Obviously, \(\fib{M}{\set{\Ouu{q}{1}}}=\nul{M}\) and \(\fib{M}{\Cq}=\Cp\).
\erem

 We write \(A^\ad\) for the conjugate transpose and \(\ran{A}\defeq\setaca{Ax}{x\in\Cq}\) for the column space of a complex \tpqa{matrix} \(A\), \tresp{}
 With the Euclidean scalar product \(\ipE{\cdot}{\cdot}\colon\x{\Cq}\to\C\) given by \(\ipE{x}{y}\defeq y^\ad x\), which is \(\C\)\nobreakdash-linear in its first argument, the vector space \(\Cq\) over the field \(\C\) becomes a unitary space.
 Let \(\mathcal{U}\) be an arbitrary non-empty subset of \(\Cq\).
 The orthogonal complement \(\mathcal{U}^\orth\defeq\setaca{v\in\Cq}{\ipE{v}{u}=0\text{ for all }u\in\mathcal{U}}\) of \(\mathcal{U}\) is a linear subspace of the unitary space \(\Cq\).
 If \(\mathcal{U}\) is a linear subspace itself, the unitary space \(\Cq\) is the orthogonal sum of \(\mathcal{U}\) and \(\mathcal{U}^\orth\).
 In this case, we write \(\OPu{\mathcal{U}}\) for the transformation matrix corresponding to the orthogonal projection onto \(\mathcal{U}\) with respect to the standard basis of \(\Cq\), \tie{}, \(\OPu{\mathcal{U}}\) is the uniquely determined matrix \(P\in\Cqq\) satisfying \(P^2=P=P^\ad\) and \(\ran{P}=\mathcal{U}\).
 For each matrix \(A\in\Cggq\), there exists a uniquely determined matrix \(Q\in\Cggq\) with \(Q^2=A\) called the \emph{\tnnH{} square root} \(Q=\sqrt{A}\) of \(A\).
 
 As a starting point of our subsequent considerations, we choose a matrix which will turn out to coincide with a matrix which was introduced in an alternate way in~\zitaa{MR2160825}{see~\ceq{5.1.11}{141} and \cthm{5.8}{}}.

\bdefnl{de1ab3}
 Let \(A\in \Cggq\) and let \(\cV \) be a linear subspace of \(\Cq\).
 Then we call the matrix \(\Psfib{A}{\cV}\defeq\OPu{\fib{\sqrt{A}}{\cV}}\) the \emph{\tPsfib{A}{\cV}} and the matrix \(\SO{A}{\cV}\defeq\sqrt{A}\Psfib{A}{\cV}\sqrt{A}\) is said to be the \emph{\tSOo{\(A\)}{\(\cV \)}}.
\edefn

 Let \(\Iq\defeq\matauuo{\Kronu{jk}}{j,k}{1}{q}\) be the identity matrix from \(\Cqq\), where \(\Kronu{jk}\) is the Kronecker delta.
 Sometimes, we will omit the indices and write \(\EM\).
 
\breml{be3ab3}
 If \(A\in\Cggq\), then \(\Psfib{A}{\set{\Ouu{q}{1}}} =\OPu{\nul{A}}\), \(\SO{A}{\set{\Ouu{q}{1}}} = \Oqq\), \(\Psfib{A}{\Cq}= \Iq \), and \(\SO{A}{\Cq} = A\).
\erem

\breml{N23}
 Let \(A\in\Cggq\) and let \(\cV\) be a linear subspace of \(\Cq\).
 Then the matrices \(\Psfib{A}{\cV}\) and \(\SO{A}{\cV}\) are both \tH{}.
\erem

\breml{be5ab3}
 If \(\cV \) is a linear subspace of \(\Cq\), then \(\Psfib{\Iq}{\cV} = \OPu{\cV}\), \(\SO{\Iq}{\cV} = \OPu{\cV}\), \(\Psfib{\OPu{\cV}}{\cV} = \Iq\), and \( \SO{\OPu{\cV}}{\cV} = \OPu{\cV}\).
\erem

 We write \(\rank A\) for the rank of a complex \tpqa{matrix} \(A\) and \(\Cgq\) for the set of \tpH{} matrices from \(\Cqq\).
 The set \(\CHq\) of \tH{} matrices from \(\Cqq\) is a partially ordered vector space over the field \(\R\) with positive cone \(\Cggq\).
 For two complex \tqqa{matrices} \(A\) and \(B\), we write \(A\lleq B\) or \(B\lgeq A\) if \(A,B\in\CHq\) and \(B-A\in\Cggq\) are fulfilled.
 For a complex \tqqa{matrix} \(A\), we have obviously \(A\lgeq\NM\) if and only if \(A\in\Cggq\).
 The above mentioned partial order \(\lleq\) on the set of \tH{} matrices is sometimes called \emph{L\"owner semi-ordering}.
 Parts of the following proposition coincide with results stated in~\zitaa{MR2160825}{\cthmss{5.3}{5.6} in combination with \cthm{5.8}{}}.
 
\bpropl{sa1ab3}
 If \(A\in\Cggq\) and \(\cV \) is a linear subspace of \(\Cq\), then:
 \benui
  \il{sa1ab3.a} \(\ran{\SO{A}{\cV}} = \ran{\sqrt{A}\Psfib{A}{\cV}}\), \(\nul{\SO{A}{\cV}} = \nul{\Psfib{A}{\cV}\sqrt{A}}\), and \(\Oqq\lleq\SO{A}{\cV}\lleq A\).
  \il{sa1ab3.b} \(\ran{\SO{A}{\cV}} = \ran{A}\cap \cV \), \(\nul{\SO{A}{\cV}} = \nul{A}+\cV^\orth  \), and in particular \( \rank{\SO{A}{\cV}} \leq \min \set{\rank{A}, \dim{\cV}}\).
  \il{sa1ab3.g} The following statements are equivalent:
  \begin{aeqii}{0}
   \il{sa1ab3.I} \(\SO{A}{\cV} = A\).
   \il{sa1ab3.II} \(\ran{\SO{A}{\cV}} = \ran{A}\).
   \il{sa1ab3.III} \(\ran{A} \subseteq \cV \).
  \end{aeqii}
  \il{sa1ab3.h} \(\ran{\SO{A}{\cV}} = \cV \) if and only if \(\cV \subseteq \ran{A}\).
  \il{sa1ab3.l} \(\SO{A}{\cV} \in \Cgq\) if and only if \(A \in \Cgq \) and \(\cV = \Cq\).
  \il{sa1ab3.k} If \(\cV \neq \Cq\), then \(\SO{A}{\cV} \in \Cggq \setminus \Cgq\).
 \eenui
\eprop
\bproof
 \eqref{sa1ab3.a} Use \(\Oqq\leq\Psfib{A}{\cV}\leq\Iq\) and \(\Psfib{A}{\cV}^2=\Psfib{A}{\cV}\).
 
 \eqref{sa1ab3.b} First we consider an arbitrary \(y\in\ran{\SO{A}{\cV}}\).
 According to~\eqref{sa1ab3.a}, we have \(y=\sqrt{A}\Psfib{A}{\cV}z\) with some \(z\in\Cq\).
 Obviously, \(x\defeq\Psfib{A}{\cV}z\) belongs to \(\fib{A}{\cV}\) and fulfills \(y=\sqrt{A}x\).
 In particular, \(y\) belongs to \(\ran{A}\cap\cV\).
 Conversely, now assume that \(y\) belongs to \(\ran{A}\cap\cV\).
 Then \(y\in\ran{\sqrt{A}}\), \tie{} there exists an \(x\in\Cq\) with \(y=\sqrt{A}x\).
 Consequently, \(\sqrt{A}x\in\cV\), \tie{} \(x\in\fib{A}{\cV}\).
 This implies \(\Psfib{A}{\cV}x=x\).
 Hence, \(y\in\ran{\sqrt{A}\Psfib{A}{\cV}}\).
 Taking~\eqref{sa1ab3.a} into account, then \(y\in\ran{\SO{A}{\cV}}\) follows.
 Thus, \(\ran{\SO{A}{\cV}} =\ran{A}\cap \cV \) is proved.
 Therefore, we get
\begin{multline*}
 \nul{\SO{A}{\cV}}
 =\ran{\SO{A}{\cV}^\ad}^\orth
 =\ran{\SO{A}{\cV}}^\orth
 =\ek*{\ran{A}\cap \cV}^\orth\\
 =\ran{A}^\orth+\cV^\orth
 =\nul{A^\ad}+\cV^\orth
 =\nul{A} +\cV^\orth.
\end{multline*}

 \eqref{sa1ab3.g} Condition~\rstat{sa1ab3.I} is obviously sufficient for~\rstat{sa1ab3.II}.
 According to~\eqref{sa1ab3.b} statements~\rstat{sa1ab3.II} and~\rstat{sa1ab3.III} are equivalent.
 If~\rstat{sa1ab3.III} holds true then \(\fib{\sqrt{A}}{\cV}=\Cq\) and, consequently, \(\Psfib{ \sqrt{A}}{\cV}=\Iq\).
 This means that~\rstat{sa1ab3.III} implies~\rstat{sa1ab3.I}.
 
 \eqref{sa1ab3.h} This equivalence follows from~\eqref{sa1ab3.b}.
 
 \eqref{sa1ab3.l} This is an immediate consequence of the definition of \(\SO{A}{\cV}\).
 
 \eqref{sa1ab3.k} This follows from~\eqref{sa1ab3.a} and~\eqref{sa1ab3.l}.
\eproof

\breml{fo1ab3}
 If \(A \in \Cgq\) and \(\cV \) is a linear subspace of \( \Cq\), then \rpropp{sa1ab3}{sa1ab3.h} shows that \(\ran{\SO{A}{\cV}} =\cV\).
\erem

 In the sequel the Moore--Penrose inverse plays an essential role.
 For this reason, we recall this notion.
 For each matrix \(A\in\Cpq\), there exists a uniquely determined matrix \(X\in\Cqp\), satisfying the four equations
\begin{align*}
 AXA&=A,&
 XAX&=X,&
 \rk{AX}^\ad&=AX,&
&\text{and}&
 \rk{XA}^\ad&=XA.
\end{align*}
 This matrix \(X\) is called the \emph{Moore--Penrose inverse of \(A\)} and is denoted by \(A^\mpi\).

 Given \(n\in\N\) and arbitrary rectangular complex matrices \(A_1,A_2,\dotsc,A_n\), we write \(\diag\seq{A_j}{j}{1}{n}=\diag\rk{A_1,A_2,\dotsc,A_n}\defeq\matauuuo{\Kronu{jk}A_j}{j}{k}{1}{n}\) for the corresponding block diagonal matrix.
 
 The following lemma yields essential insights into the structure of the Schur complement with respect to the \tSOo{\(A\)}{\(\cV \)}.
 
\bleml{lm2ab3}
 Assume \(q\geq2\), let \(d \in \mn{1}{ q-1 } \), and let \(\cV \) be a linear subspace of \(\Cq\) with \( \dim\cV = d\).
 Let \(u_1,u_2,\dotsc,u_q\) be an orthonormal basis of \(\Cq\) such that \(u_1,u_2,\dotsc,u_d\) is an orthonormal basis of \(\cV \) and let \( U \defeq\mat{u_{1}, \dotsc, u_{q}}\).
 Let \(A \in \Cggq\) and let \(B=\matauuo{B_{jk}}{j,k}{1}{2}\) be the \tbr{} of \(B \defeq U^\ad A U \) with \taaa{d}{d}{block} \(B_{11}\).
 Then
\beql{lm2ab3.A}
 U^\ad \SO{A}{\cV} U
 = \diag\rk{B_{11} - B_{12}B_{22}^\mpi B_{21}, \Ouu{(q-d)}{(q-d)}}.
\eeq
\elem
\bproof
 Let \(R\defeq U^\ad\sqrt{A}U\) and let \(R=\tmatp{F}{G}\) be the \tbr{} of \(R\) with \taaa{d}{q}{block} \(F\).
 Then
\(
 B
 =RR^\ad
 =
 \tmat{
  EE^\ad&EF^\ad\\
  FE^\ad& FF^\ad
 }
\).
 Consequently,
\beql{lm2ab3.1}
 B_{11} - B_{12}B_{22}^\mpi B_{21}
 =E\ek*{\Iq-F^\ad\rk{FF^\ad}^\mpi F}E^\ad
 =E\rk{\Iq-F^\mpi F}E^\ad.
\eeq
 We set \(U_1\defeq\mat{u_1,\dotsc,u_d}\)  and \(\cW\defeq\fib{\sqrt{A}}{\cV}\).
 Because of \(U^\ad U=\Iq\), we have
\[\begin{split}
  U^\ad\cW 
  &=\setaca{U^\ad x}{x\in\Cq\text{ and }\sqrt{A}x\in\cV}
  =\setaca{y\in\Cq}{\sqrt{A}Uy\in\cV}\\
  &=\setaca*{y\in\Cq}{\sqrt{A}Uy\in\ran{U_1}}
  =\setava*{y\in\Cq}{\exists z\in\Co{d}\colon\sqrt{A}Uy=U\tmatp{z}{\Ouu{\rk{q-d}}{1}}}\\
  &=\setaca*{y\in\Cq}{\mat{\Ouu{\rk{q-d}}{d},\Iu{q-d}}U^\ad\sqrt{A}Uy=\Ouu{\rk{q-d}}{1}}
  =\nul{F}.
\end{split}\]
 Hence, \(\OPu{U^\ad\cW }=\Iq-F^\mpi F\).
 Since \(U^\ad\OPu{\cW }U\) is an idempotent and \tH{} complex matrix fulfilling \(\ran{U^\ad\OPu{\cW }U}=U^\ad\cW \), we get then \(U^\ad\OPu{\cW }U=\Iq-F^\mpi F\).
 Thus, in view of \(\OPu{\cW}=\Psfib{A}{\cV}\), we obtain
\begin{multline*}
 U^\ad\SO{A}{\cV}U
 =\rk{U^\ad\sqrt{A}U}\rk{U^\ad\OPu{\cW}U}\rk{U^\ad\sqrt{A}U}^\ad\\
 =R\rk{\Iq-F^\mpi F}R^\ad
 =
 \begin{pmat}[{|}]
  E\rk{\Iq-F^\mpi F}E^\ad&E\rk{\Iq-F^\mpi F}F^\ad\cr\-
  F\rk{\Iq-F^\mpi F}E^\ad&F\rk{\Iq-F^\mpi F}F^\ad\cr
 \end{pmat}.
\end{multline*}
 Using \eqref{lm2ab3.1}, \(F\rk{\Iq-F^\mpi F}=\NM\), and \(F^\mpi F=\rk{F^\mpi F}^\ad\), then \eqref{lm2ab3.A} follows.
\eproof

 The following observation makes clear why in \rdefn{de1ab3} the terminology ``\tSOo{\(A\)}{\(\cV \)}'' was chosen.

\breml{be7ab3}
 Assume \(q\geq2\), let \(d \in \mn{1}{q-1}\), let \(A=\matauuo{A_{jk}}{j,k}{1}{2}\) be the \tbr{} of a matrix \(A \in \Cggq\) with \taaa{d}{d}{block} \(A_{11}\), and let \(\cV\defeq\ran{\tmatp{\Iu{d}}{\Ouu{(q-d)}{d}}}\).
 Then
\(
 \SO{A}{\cV}= \diag\rk{A_{11} - A_{12}A_{22}^\mpi A_{21}, \Ouu{(q-d)}{(q-d)}}
\).
\erem

 The following result shows in combination with~\zitaa{MR2160825}{\cthm{5.1}{}} that the construction introduced in \rdefn{de1ab3} coincides with the matrix introduced by Ando~\zitaa{And05}{\cfo{5.1.11}{141}}.
 
\bpropl{sa2ab3}
 Let \(A \in \Cggq\) and let \(\cV \) be a linear subspace of \(\Cq\).
 For all \(x\in\Cq\), then
\beql{sa2ab3.A}
  x^\ad\SO{A}{\cV} x 
  = \min_{y \in\cV^\orth} \rk{x-y}^\ad A \rk{x-y}.
\eeq
\eprop
\bproof
 Set \(d\defeq\dim\cV\).
 We consider an arbitrary \(x\in\Cq\).
 If \(d=0\), then \rpropp{sa1ab3}{sa1ab3.b} yields \(\SO{A}{\cV}=\Oqq\) and, because of \(A\in\Cggq\), therefore \eqref{sa2ab3.A}.
 If \(d=q\), then \(\cV^\orth=\set{\NM}\) and \rpropp{sa1ab3}{sa1ab3.g} yields \(\SO{A}{\cV}=A\), which implies \eqref{sa2ab3.A}.
 Now we consider the case \(1\leq d\leq q-1\).
 Let \(u_1,u_2,\dotsc,u_q\) be an orthonormal basis of \(\Cq\) such that \(u_1,u_2,\dotsc,u_d\) is an orthonormal basis of \(\cV \).
 Let \(U_1\defeq\mat{u_{1}, \dotsc, u_d}\) and let \(U_2\defeq\mat{u_{d+1}, \dotsc, u_q}\).
 Then \(U\defeq\mat{U_1,U_2}\) is unitary.
 Let \(B=\matauuo{B_{jk}}{j,k}{1}{2}\) be the \tbr{} of \(B \defeq U^\ad A U \) with \taaa{d}{d}{block} \(B_{11}\).
 Setting \(S\defeq B_{11} - B_{12}B_{22}^\mpi B_{21}\), from \rlem{lm2ab3.A} we get \(U^\ad \SO{A}{\cV} U = \diag\rk{S, \Ouu{(q-d)}{(q-d)}}\).
 Let \(f_1\defeq U_1^\ad x\) and let \(f_2\defeq U_2^\ad x\).
 Then \(f\defeq U^\ad x\) admits the \tbr{} \(f=\tmatp{f_1}{f_2}\).
 Thus, we get
\beql{sa2ab3.1}
 x^\ad\SO{A}{\cV}x
 =f^\ad\ek*{\diag\rk{S, \Ouu{(q-d)}{(q-d)}}}f
 =f_1^\ad Sf_1.
\eeq
 We consider an arbitrary \(y\in\cV^\orth\).
 Setting \(g_1\defeq U_1^\ad y\) and \(g_2\defeq U_2^\ad y\), we see that \(g\defeq U^\ad y\) admits the \tbr{} \(g=\tmatp{g_1}{g_2}\) and that \(g_1=\Ouu{d}{1}\).
 Thus, \(h\defeq f-g\) can be represented via \(h=\tmatp{h_1}{h_2}\), where \(h_1\defeq f_1\) and \(h_2\defeq f_2-g_2\).
 The matrices \(B\) and \(B_{22}\) are obviously \tnnH{}.
 Therefore, using a well-known factorization formula (see, \teg{}~\zitaa{MR1152328}{\clemss{1.1.9}{18}{1.1.7}{17}}), we have \(B=R^\ad\ek{\diag\rk{S,B_{22}}}R\) where \(R\defeq\smat{\Iu{d}&\Ouu{d}{\rk{q-d}}\\B_{22}^\mpi B_{21}&\Iu{q-d}}\).
 It is easily checked that \(Rh=\tmatp{f_1}{r_2-g_2}\) where \(r_2\defeq B_{22}^\mpi B_{21}f_1+f_2\).
 Applying \(x-y=Uh\) and \eqref{sa2ab3.1} we conclude then
\begin{multline*}
 \rk{x-y}^\ad A\rk{x-y}
 =h^\ad Bh
 =\matp{f_1}{r_2-g_2}^\ad\ek*{\diag\rk{S,B_{22}}}\matp{f_1}{r_2-g_2}\\
 =f_1^\ad Sf_1-\rk{r_2-g_2}^\ad B_{22}\rk{r_2-g_2}
 =x^\ad\SO{A}{\cV}x-\rk{r_2-g_2}^\ad B_{22}\rk{r_2-g_2}.
\end{multline*}
 Consequently, \(\rk{x-y}^\ad A\rk{x-y}\geq x^\ad\SO{A}{\cV}x\) for all \(y\in\cV^\orth\) with equality if and only if \(r_2-g_2\in\nul{B_{22}}\).
 Obviously, the particular vector \(y\defeq U_2r_2\) belongs to \(\cV^\orth\) and we have \(g_2=r_2\).
 Thus, \(y\defeq U_2r_2\) fulfills \(B_{22}\rk{r_2-g_2}=\NM\).
\eproof

 \rprop{sa2ab3} establishes the coincidence of the matrix introduced in \rdefn{de1ab3} with the construction used by T.~Ando~\zitaa{And05}{\cfo{5.1.11}{141}}.
 This leads us to interesting insights about the main object of this section.

  For two \tH{} \tqqa{matrices} \(A\) and \(B\) with \(A\lleq B\), the (closed) matricial interval \(\matint{A}{B}\defeq\setaca{X\in\CHq}{A\lleq X\lleq B}\) is non-empty.
  
\bnotal{N1452}
 If \(A \in \Cggq\) and \(\cV \) is a linear subspace of \(\Cq\), then let \(\LcR{A}{\cV}\defeq\setaca{X\in\matint{\Oqq}{A}}{\ran{X}\subseteq \cV}\).
\enota

\bthmnl{\tcf{}~\zitaa{MR2160825}{\cthm{5.3}{}}}{sa3ab3}
 If \(A \in \Cggq\) and \(\cV \) is a linear subspace of \(\Cq\), then \(\SO{A}{\cV} \in \LcR{A}{\cV} \) and \(\SO{A}{\cV}\lgeq X\) for all \(X\in \LcR{A}{\cV} \).
\ethm
\bproof
 For the convenience of the reader, we reproduce the proof given in~\zitaa{MR2160825}{\cthm{5.3}{}}:
 From \rpartss{sa1ab3.a}{sa1ab3.b} of \rprop{sa1ab3} we infer \(\SO{A}{\cV}\in\LcR{A}{\cV}\).
 Now consider an arbitrary \(X\in\LcR{A}{\cV}\).
 Furthermore, consider an arbitrary \(x\in\Cq\) and an arbitrary \(y_0\in\cV^\orth\).
 Because of \(\ran{X}\subseteq\cV\) and \(X\in\Cggq\), then \(y_0\in\ran{X}^\orth=\nul{X^\ad}=\nul{X}\), implying \(\rk{x-y_0}^\ad X\rk{x-y_0}=x^\ad Xx\).
 Because of \(A\lgeq X\), moreover \(\rk{x-y_0}^\ad A\rk{x-y_0}\lgeq\rk{x-y_0}^\ad X\rk{x-y_0}\).
 Taking additionally into account \rprop{sa2ab3}, we thus obtain \eqref{sa2ab3.A} and, hence,
\(
 x^\ad\SO{A}{\cV}x
 \lgeq x^\ad Xx
\).
 Consequently, \(\SO{A}{\cV}\lgeq X\).
\eproof

\breml{A.R1004}
 Let \(A,B\in\Cggq\).
 Then \(\ran{B}\subseteq\ran{A}\) if and only if \(B\lleq\gamma A\) for some \(\gamma\in(0,\infp)\).
 In this case, \(\gamma_0\defeq\sup\setaca{\rk{x^\ad Bx}\rk{x^\ad A x}^\inv}{x\in\Cq\setminus\nul{A}}\) fulfills \(\gamma_0\in[0,\infp)\) and \(B\lleq\gamma_0A\).
\erem

\blemnl{\tcf{}~\zitaa{MR2160825}{Equivalence~(5.0.7)}}{A.L0925}
 If \(A,B\in\Cggq\), then \(\ran{A}\cap\ran{B}=\set{\Ouu{q}{1}}\) if and only if \(\matint{\Oqq}{A}\cap\matint{\Oqq}{B}=\set{\Oqq}\).
\elem
\bproof
 Let \(A,B\in\Cggq\).
 Then \(\Oqq\in\matint{\Oqq}{A}\cap\matint{\Oqq}{B}\).
 
 First assume \(\ran{A}\cap\ran{B}=\set{\Ouu{q}{1}}\).
 We consider an arbitrary \(X\in\matint{\Oqq}{A}\cap\matint{\Oqq}{B}\).
 Then \(\ran{X}\subseteq\ran{A}\cap\ran{B}\), by virtue of \rrem{A.R1004}.
 Consequently, \(X=\Oqq\).
 
 Conversely, assume that \(\ran{A}\cap\ran{B}\neq\set{\Ouu{q}{1}}\).
 Then \(P\defeq\OPu{\ran{A}\cap\ran{B}}\) fulfills \(P\in\Cggq\setminus\set{\Oqq}\) and \(\ran{P}=\ran{A}\cap\ran{B}\).
 From \rrem{A.R1004} we can conclude then the existence \(\alpha,\beta\in(0,\infp)\) with \(P\lleq\alpha A\) and \(P\lleq\beta B\).
 Then \(\gamma\defeq\min\set{1/\alpha,1/\beta}\) fulfills \(\gamma\in(0,\infp)\) and \(\gamma P\in\rk{\matint{\Oqq}{A}\cap\matint{\Oqq}{B}}\setminus\set{\Oqq}\).
\eproof

\bthmnl{\tcf{}~\zitaa{MR2160825}{\cthm{5.7}{}}}{A.T1514}
 Let \(A \in \Cggq\) and let \(\cV \) be a linear subspace of \(\Cq\).
 Then:
\benui
 \il{A.T1514.a} \(\ran{A-\SO{A}{\cV}}\cap\cV=\set{\Ouu{q}{1}}\).
 \il{A.T1514.b} Let \(X,Y\in\Cggq\) be such that \(X+Y=A\).
 Then the following statements are equivalent:
\begin{aeqii}{0}
 \il{A.T1514.i} \(\ran{X}\subseteq\cV\) and \(\ran{Y}\cap\cV=\set{\Ouu{q}{1}}\).
 \il{A.T1514.ii} \(X=\SO{A}{\cV}\) and \(Y=A-\SO{A}{\cV}\).
\end{aeqii}
\eenui
\ethm
\bproof
 For the convenience of the reader, we reproduce the proof given in~\zitaa{MR2160825}{\cthm{5.7}{}}:
 First observe that \(\Oqq\lleq\SO{A}{\cV}\lleq A\) and \(\ran{\SO{A}{\cV}}\subseteq\cV\), by virtue of \rpartss{sa1ab3.a}{sa1ab3.b} of \rprop{sa1ab3}.
 
 \eqref{A.T1514.a} We have \(\set{A-\SO{A}{\cV},\OPu{\cV}}\subseteq\Cggq\).
 Hence, \(\Oqq\in\matint{\Oqq}{A-\SO{A}{\cV}}\cap\matint{\Oqq}{\OPu{\cV}}\).
 Consider an arbitrary \(X\in\matint{\Oqq}{A-\SO{A}{\cV}}\cap\matint{\Oqq}{\OPu{\cV}}\).
 Then \(X\in\Cggq\) and, moreover, \(X\lleq A-\SO{A}{\cV}\) and \(X\lleq\OPu{\cV}\).
 In particular, \(\ran{X}\subseteq\ran{\OPu{\cV}}=\cV\), by virtue of \rrem{A.R1004}.
 Taking additionally into account \(\SO{A}{\cV}\in\Cggq\) and \(\ran{\SO{A}{\cV}}\subseteq\cV\), we thus obtain \(\SO{A}{\cV}+X\in\LcR{A}{\cV}\).
 Hence, \rthm{sa3ab3} yields \(\Oqq\lleq\SO{A}{\cV}+X\lleq\SO{A}{\cV}\).
 This implies \(X=\Oqq\).
 Consequently, \(\matint{\Oqq}{A-\SO{A}{\cV}}\cap\matint{\Oqq}{\OPu{\cV}}=\set{\Oqq}\).
 According to \rlem{A.L0925}, then \(\ran{A-\SO{A}{\cV}}\cap\ran{\OPu{\cV}}=\set{\Ouu{q}{1}}\) follows.
 In view of \(\ran{\OPu{\cV}}=\cV\), the proof of~\rpart{A.T1514.a} is complete.

 \eqref{A.T1514.b} Obviously \(\Oqq\lleq X\lleq A\).
 
 First suppose~\rstat{A.T1514.i}.
 Then \(X\in\LcR{A}{\cV}\).
 \rthm{sa3ab3} yields then \(\SO{A}{\cV}\lgeq X\).
 Taking into account \(A-\SO{A}{\cV}\in\Cggq\), thus \(\Oqq\lleq\SO{A}{\cV}-X\lleq A-X\).
 According to \rrem{A.R1004} and \(Y=A-X\), then \(\ran{\SO{A}{\cV}-X}\subseteq\ran{Y}\).
 Because of \(X\in\Cggq\), we can furthermore conclude \(\Oqq\lleq\SO{A}{\cV}-X\lleq\SO{A}{\cV}\).
 By virtue of \rrem{A.R1004} and \(\ran{\SO{A}{\cV}}\subseteq\cV\), then \(\ran{\SO{A}{\cV}-X}\subseteq\ran{\SO{A}{\cV}}\subseteq\cV\).
 Taking additionally into account~\rstat{A.T1514.i}, we thus obtain \(\ran{\SO{A}{\cV}-X}\subseteq\ran{Y}\cap\cV=\set{\Ouu{q}{1}}\).
 Consequently, \(\SO{A}{\cV}=X\).
 Thus,~\rstat{A.T1514.ii} holds true.

 If we conversely suppose~\rstat{A.T1514.ii}, then~\rstat{A.T1514.i} follows from \(\ran{\SO{A}{\cV}}\subseteq\cV\) and~\eqref{A.T1514.a}.
\eproof

\section{On a restricted extension problem for a finite \tHnnde{} sequence}\label{S.H}
 This section is written against to the background of applying the results of \rsec{S.ando} to the truncated matricial Hamburger moment problems formulated in the introduction for \(\Omega=\R\).
 In the heart of our strategy lies the treatment of a special restricted extension problem for matrices.
 The complete answer to this problem is contained in \rthm{sa4ab5} which is the central result of this section.

 If \(\kappa\in\NOinf\) and \(\seqska\) is a sequence of complex \tpqa{matrices}, then let \(\Hu{n}\defeq\matauuuo{s_{j+k}}{j}{k}{0}{n}\) for all \(n\in\NO\) with \(2n\leq\kappa\).

 For each \(n\in\NO\), by \(\Hggq{2n}\) we denote the set of all sequences \(\seqs{2n}\) of complex \tqqa{matrices} for which the corresponding \tbHm{} \(\Hu{n}\)
 is \tnnH{}.
 Furthermore, denote by \(\Hggqinf\) the set of all sequences \(\seqsinf\) of complex \tqqa{matrices} satisfying \(\seqs{2n}\in\Hggq{2n}\) for all \(n\in\NO\).
 The sequences belonging to \(\Hggq{2n}\) or \(\Hggqinf\) are said to be \emph{\tHnnd}.
 Using~\zitaa{MR3014199}{\clem{3.2}{126}}, we can conclude:
 
\breml{H.R0857}
 If \(\kappa\in\NOinf\) and if \(\seqs{2\kappa}\in\Hggq{2\kappa}\), then \(\su{j}\in\CHq\) for all \(j\in\mn{0}{2\kappa}\) and \(\su{2k}\in\Cggq\) for all \(k\in\mn{0}{\kappa}\).
\erem

 Given \(n\in\N\) arbitrary rectangular complex matrices \(A_1,A_2,\dotsc,A_n\), we write \(\col\seq{A_j}{j}{1}{n}=\col\rk{A_1,A_2,\dotsc,A_n}\) (\tresp{}, \(\row\seq{A_j}{j}{1}{n}\defeq\mat{A_1,A_2,\dotsc,A_n}\)) for the block column (\tresp{}, block row) build from the matrices \(A_1,A_2,\dotsc,A_n\) if their numbers of columns (\tresp{}, rows) are all equal.

\bnotal{H.N.yz}
 Let \(\kappa\in\NOinf\) and let \(\seqska \) be a sequence of complex \tpqa{matrices}.
\benui
 \il{H.N.yz.a} Let \(\yuu{\ell}{m}\defeq\col\seq{s_j}{j}{\ell}{m}\) and \(\zuu{\ell}{m}\defeq\row\seq{s_j}{j}{\ell}{m}\)
 for all \(\ell,m\in\NO\) with \(\ell\leq m\leq\kappa\).
 \il{H.N.yz.b} Let \(\Tripu{0}\defeq\Opq\) and let \(\Tripu{n}\defeq\zuu{n}{2n-1}\Hu{n-1}^\mpi\yuu{n}{2n-1}\) for all \(n\in\N\) with \(2n-1\leq\kappa\).
 \il{H.N.yz.c} Let \(\Lu{n}\defeq\su{2n}-\Tripu{n}\) for all \(n\in\NO\) with \(2n\leq\kappa\).
\eenui
\enota

\breml{H.R.L>}
 If \(\kappa\in\NOinf\) and \(\seqs{2\kappa}\in\Hggq{2\kappa}\), then~\zitaa{MR2570113}{\crem{2.1(a)}{765}} shows that \(\seqs{2n}\in\Hggq{2n}\) and \(\Lu{n}\in\Cggq\) for all \(n\in\mn{0}{\kappa}\).
\erem

 Let \(n\in\NO\).
 Denote by \(\Hggeq{2n}\) the set of all sequences \(\seqs{2n}\) of complex \tqqa{matrices} for which there exists a pair \((\su{2n+1},\su{2n+2})\) of complex \tqqa{matrices} such that the sequence \(\seqs{2n+2}\) belongs to \(\Hggq{2n+2}\).
 Denote by \(\Hggeq{2n+1}\) the set of all sequences \(\seqs{2n+1}\) of complex \tqqa{matrices} for which there exists a complex \tqqa{matrix} \(\su{2n+2}\) such that the sequence \(\seqs{2n+2}\) belongs to \(\Hggq{2n+2}\).
 Furthermore, let \(\Hggeqinf\defeq\Hggqinf\).
 The sequences belonging to \(\Hggeq{2n}\), \(\Hggeq{2n+1}\), or \(\Hggeqinf\) are said to be \emph{\tHnnde}.
 
\breml{H.R.He<gg}
 \rrem{H.R.L>} shows that \(\Hggeq{2\kappa}\subseteq\Hggq{2\kappa}\) for all \(\kappa\in\NOinf\).
\erem
 
\breml{H.R0904}
 If \(\kappa\in\NOinf\) and \(\seqska\in\Hggeq{\kappa}\), then, because of \rremss{H.R.L>}{H.R0857}, we have \(\Tripu{n}\in\Cggq\) for all \(n\in\mn{0}{\kappa}\).
\erem

\bnotal{H.N1949}
 Let \(n\in\N\) and let \(\seqs{2n-1}\) is a sequence of complex \tqqa{matrices}.
 Then let \(\exHggs{2n-1}\) be the set of all \(\su{2n}\in\Cqq\) for which \(\seqs{2n}\) belongs to \(\Hggq{2n}\).
 Obviously, \(\exHggs{2n-1}\neq\emptyset\) if and only if \(\seqs{2n-1}\in\Hggeq{2n-1}\).
\enota

\bpropnl{\tcf{}~\zitaa{MR2570113}{\cprop{2.22(a)}{772}}}{H.P0840}
 Let \(n\in\N\), let \(\seqs{2n-1}\in\Hggeq{2n-1}\), and let \(\su{2n}\in\Cqq\).
 Then  \(\su{2n}\in\exHggs{2n-1}\) if and only if \(\Lu{n}\in\Cggq\).
\eprop

\bnotal{H.N0843}
 Let \(n\in\N\), let \(\seqs{2n-1}\) be a sequence of complex \tqqa{matrices}, and let \(Y\in\CHq\).
 Then denote by \(\exHggsk{2n-1}{Y}\) the set of all \(X\in\exHggs{2n-1}\) satisfying \(Y-X\in\Cggq\).
\enota

\bleml{H.L0845}
 Let \(\kappa\in\NOinf\) and let \(\seqs{2\kappa}\in\Hggq{2\kappa}\).
 Then \(\Oqq\lleq\Tripu{n}\lleq\su{2n}\) for all \(n\in\mn{0}{\kappa}\).
\elem
\bproof
 \rrem{H.R0857} shows that \(\Oqq\lleq\Tripu{0}\lleq\su{0}\).
 Now assume \(\kappa\geq1\) and consider an arbitrary \(n\in\mn{1}{\kappa}\).
 \rrem{H.R0857} yields \(\su{2n}\in\CHq\).
 According to \rrem{H.R.L>}, furthermore \(\Lu{n}\in\Cggq\) and \(\seqs{2n}\in\Hggq{2n}\).
 In particular, \(\seqs{2n-1}\in\Hggeq{2n-1}\).
 Thus, \rrem{H.R0904} and \rnotap{H.N.yz}{H.N.yz.c} yield \(\Oqq\lleq\Tripu{n}\lleq\su{2n}\).
\eproof

 To indicate that a certain (block) matrix \(X\) is built from a sequence \(\seqska\), we sometimes write \(X^{\rk{s}}\) for \(X\).
 
\bpropl{H.P0850}
 If \(n\in\N\) and \(\seqs{2n}\in\Hggq{2n}\), then \(\Oqq\lleq\Tripu{n}\lleq\su{2n}\) and \(\exHggsk{2n-1}{\su{2n}}=\matint{\Tripu{n}}{\su{2n}}\).
\eprop
\bproof
 By virtue of \rlem{H.L0845}, we have \(\Oqq\lleq\Tripu{n}\lleq\su{2n}\).
 In particular, the matrices \(\Tripu{n}\) and \(\su{2n}\) are \tH{} .
 Let the sequence \(\seq{t_j}{j}{0}{2n-1}\) be given by \(t_j\defeq\su{j}\) for each \(j\in\mn{0}{2n-1}\).
 Obviously, then \(\Tripuo{n}{t}=\Tripu{n}\) and \(\set{\seqs{2n-1},\seq{t_j}{j}{0}{2n-1}}\subseteq\Hggeq{2n-1}\).
 Obviously, then
 \beql{M11}
 t_{2n}-\Tripu{n}
 =t_{2n}-\Tripuo{n}{t}
 =\Luo{n}{t}.
\eeq

 We first consider an arbitrary \(t_{2n}\in\exHggsk{2n-1}{\su{2n}}\).
 Then \(\seqt{2n}\in\Hggq{2n}\) and \(\su{2n}-t_{2n}\in\Cggq\).
 In particular, \(t_{2n}\in\CHq\).
 Since \rrem{H.R.L>} shows that \(\Luo{n}{t}\) is \tnnH{}, \eqref{M11} yields \(\Tripu{n}\lleq t_{2n}\lleq\su{n}\).
 
 Conversely, we consider now an arbitrary \(t_{2n}\in\matint{\Tripu{n}}{\su{n}}\).
 Then \(\su{2n}-t_{2n}\in\Cggq\) and \(t_{2n}-\Tripu{n}\in\Cggq\).
 Thus, \eqref{M11} yields \(\Luo{n}{t}\in\Cggq\).
 According to \rprop{H.P0840}, thus \(\seqt{2n}\in\Hggq{2n}\).
 Consequently, \(t_{2n}\in\exHggsk{2n-1}{\su{2n}}\).
\eproof

\bcorl{H.C1530}
 If \(n\in\NO\) and \(\seqs{2n}\in\Hggq{2n}\), then \(\setaca{\mpm{\sigma}{2n}}{\sigma\in\MggqRskg{2n}}\subseteq\matint{\Tripu{n}}{\su{2n}}\).
\ecor
\bproof
 In view of \(\Tripu{0}=\Oqq\), the case \(n=0\) is obvious.
 If \(n\geq1\), combine \rthm{T1708}, \rrem{H.R.He<gg}, and \rprop{H.P0850}.
\eproof

\bnotal{bz1ab5}
 If \(n\in\N\) and \(\seqs{2n-1}\) is a sequence of complex \tqqa{matrices}.
 Then let \(\exHgges{2n-1}\) be the set of all complex \tqqa{matrices} \(\su{2n}\) such that \(\seqs{2n}\in\Hggeq{2n}\).
\enota

 From~\zitaa{MR2570113}{\cprop{2.22(a), (b)}{772}} we know that \(\exHgges{2n-1}\neq\emptyset\) if and only if \(\seqs{2n-1}\) belongs to \(\Hggeq{2n-1}\).

\bpropnl{\tcf{}~\zitaa{MR2570113}{\cprop{2.22(b)}{772}}}{sa2ab5}
 Let \(n\in\N\), let \(\seqs{2n-1}\in\Hggeq{2n-1}\), and let \(\su{2n}\in\Cqq\).
 Then  \(\su{2n}\in\exHgges{2n-1}\) if and only if \(\Lu{n}\in\Cggq\) and 
\(\ran{\Lu{n}}\subseteq\ran{\Lu{n-1}}\).
\eprop

\bnotal{N1418}
 Let \(n\in\N\), let \(\seqs{2n-1}\) be a sequence of complex \tqqa{matrices}, and let \(Y\in\CHq\).
 Then denote by \(\exHggesk{2n-1}{Y}\) the set of all \(X\in\exHgges{2n-1}\) satisfying \(Y-X\in\Cggq\).
\enota

 Observe that the following construction is well defined, due to \rrem{H.R.L>} and \rdefn{de1ab3}:
 
\bnotal{H.N.F}
 Let \(\kappa\in\NOinf\) and let \(\seqs{2\kappa}\in\Hggq{2\kappa}\).
 Then let \(\RLu{0}\defeq\su{0}\).
 If \(n\in\mn{1}{\kappa}\), then let \(\RLu{n}\defeq\Tripu{n}+\SO{\Lu{n}}{\ran{\Lu{n-1}}}\).
\enota

\bleml{H.L0901}
 Let \(\kappa\in\NOinf\) and let \(\seqs{2\kappa}\in\Hggq{2\kappa}\).
 For all \(n\in\mn{0}{\kappa}\), then
\beql{M3}
 \Oqq
 \lleq\Tripu{n}
 \lleq\RLu{n}
 \lleq\su{2n}.
\eeq
\elem
\bproof
 According to \rlem{H.L0845}, we have \(\Oqq\lleq\Tripu{n}\lleq\su{2n}\) for all \(n\in\mn{0}{\kappa}\).
 In particular, the matrices \(\Tripu{n}\) and \(\su{2n}\) are \tH{} for all \(n\in\mn{0}{\kappa}\). 
 By virtue of \(\RLu{0}=\su{0}\), we get \eqref{M3} for \(n=0\).
 Now assume \(\kappa\geq1\) and \(n\in\mn{1}{\kappa}\).
 According to \rrem{H.R.L>}, we have \(\Lu{n}\in\Cggq\).
 From \rpropp{sa1ab3}{sa1ab3.a} we obtain \(\Oqq\lleq\SO{\Lu{n}}{\ran{\Lu{n-1}}}\lleq\Lu{n}\), implying, by virtue of \rnota{H.N.F} and \rnotap{H.N.yz}{H.N.yz.c} then \(\Oqq\lleq\Tripu{n}\lleq\RLu{n}\lleq\su{2n}\).
\eproof

\bthml{sa4ab5}
 Let \(n\in\N\) and let \(\seqs{2n}\in\Hggq{2n}\).
 Then:
 \benui
  \il{sa4ab5.a} \eqref{M3} and \(\exHggesk{2n-1}{\su{2n}}=\matint{\Tripu{n}}{\RLu{n}}\) hold true.
  \il{sa4ab5.b} \(\RLu{n}=\su{2n}\) if and only if \(\seqs{2n}\in\Hggeq{2n}\).
 \eenui
\ethm
\bproof
 \eqref{sa4ab5.a} From \rlem{H.L0901} we get \eqref{M3}.
 In particular, the matrices \(\Tripu{n}\), \(\RLu{n}\), and \(\su{2n}\) are \tH{}.
 \rrem{H.R.L>} yields \(\Lu{n}\in\Cggq\).
 With \(A\defeq\Lu{n}\) and \(\cV\defeq\ran{\Lu{n-1}}\) we have \(\RLu{n}=\Tripu{n}+\SO{A}{\cV}\) and \(\su{2n}=\Tripu{n}+A\).
 Let the sequence \(\seq{t_j}{j}{0}{2n-1}\) be given by \(t_j\defeq\su{j}\) for all \(j\in\mn{0}{2n-1}\).
 Then \(\Luo{n-1}{t}=\Lu{n-1}\) and \(\Tripuo{n}{t}=\Tripu{n}\).
 In particular, \(\cV=\ran{\Luo{n-1}{t}}\).
 Furthermore, \(\set{\seqs{2n-1},\seq{t_j}{j}{0}{2n-1}}\subseteq\Hggeq{2n-1}\).
 
 First consider now an arbitrary \(t_{2n}\in\exHggesk{2n-1}{\su{2n}}\).
 Then \(\seqt{2n}\in\Hggeq{2n}\) and \(\su{2n}-t_{2n}\in\Cggq\).
 In particular, \(t_{2n}\in\CHq\).
 Let \(X\defeq\Luo{n}{t}\).
 According to \rprop{sa2ab5}, then  \(X\in\Cggq\) and \(\ran{X}\subseteq\cV\).
 In view of \rnotap{H.N.yz}{H.N.yz.c}, we have \(X=t_{2n}-\Tripuo{n}{t}=t_{2n}-\Tripu{n}\).
 Hence, \(X\leq\su{2n}-\Tripu{n}=A\).
 Consequently, \(X\in\LcR{A}{\cV}\).
 \rthm{sa3ab3} yields then \(\SO{A}{\cV}\lgeq X\).
 In view of \rnotap{H.N.yz}{H.N.yz.c}, thus \(\Tripu{n}\lleq t_{2n}\lleq\RLu{n}\).
 
 Conversely, let \(t_{2n}\in\matint{\Tripu{n}}{\RLu{n}}\).
 Then \(X\defeq t_{2n}-\Tripu{n}\) is \tH{} and fulfills \(\Oqq\lleq X\lleq\SO{A}{\cV}\).
 Hence, \(\ran{X}\subseteq\ran{\SO{A}{\cV}}\) and, according to \rprop{sa1ab3}, furthermore \(\ran{\SO{A}{\cV}}=\ran{A}\cap\cV\) and \(\SO{A}{\cV}\lleq A\).
 Consequently, \(\ran{X}\subseteq\cV\) and \(\Oqq\lleq X\lleq A\).
 In view of \rnotap{H.N.yz}{H.N.yz.c}, furthermore \(\Luo{n}{t}=t_{2n}-\Tripuo{n}{t}=t_{2n}-\Tripu{n}=X\).
 Therefore, \(\Luo{n}{t}\in\Cggq\) and \(\ran{\Luo{n}{t}}\subseteq\ran{\Luo{n-1}{t}}\).
 According to \rprop{sa2ab5}, thus \(\seqt{2n}\in\Hggeq{2n}\).
 Moreover, because of \(\su{2n}-t_{2n}=A-X\), we get \(t_{2n}\lleq\su{2n}\).
 Consequently, \(t_{2n}\in\exHggesk{2n-1}{\su{2n}}\).
 
 \eqref{sa4ab5.b} In view of \(\seqs{2n}\in\Hggq{2n}\), we have \(\seqs{2n-1}\in\Hggeq{2n-1}\) and \(\su{2n}\in\exHggs{2n-1}\).
 \rprop{H.P0840} yields then \(\Lu{n}\in\Cggq\).
 By virtue of \rnota{H.N.F} and \rnotap{H.N.yz}{H.N.yz.c} we have \(\RLu{n}=\su{2n}\) if and only if \(\Lu{n}=\SO{\Lu{n}}{\ran{\Lu{n-1}}}\).
 According to \rpropp{sa1ab3}{sa1ab3.g}, the latter is equivalent to \(\ran{\Lu{n}}\subseteq\ran{\Lu{n-1}}\).
 In view of \(\seqs{2n-1}\in\Hggeq{2n-1}\) and \(\Lu{n}\in\Cggq\), \rprop{sa2ab5} shows the equivalence of \(\ran{\Lu{n}}\subseteq\ran{\Lu{n-1}}\) and \(\su{2n}\in\exHgges{2n-1}\).
 Hence,~\eqref{sa4ab5.b} follows.
\eproof

 \rthm{sa4ab5} leads us now quickly in an alternative way to one of the main results of~\zita{MR2570113}.

\bthmnl{\tcf{}~\zitaa{MR2570113}{\cthm{7.8}{809}}}{H.C1524}
 If \(n\in\NO\) and \(\seqs{2n}\in\Hggq{2n}\), then \(\setaca{\mpm{\sigma}{2n}}{\sigma\in\MggqRskg{2n}}=\matint{\Tripu{n}}{\RLu{n}}\).
\ethm
\bproof
 In view of \(\Tripu{0}=\Oqq\) and \(\RLu{0}=\su{0}\), the case \(n=0\) is obvious.
 If \(n\geq1\), combine \rthmss{T1708}{sa4ab5}.
\eproof

\section{On equivalence classes of truncated matricial moment problems of type \mprob{\R}{2n}{\lleq}}\label{H.Seq}
 This section contains an aspect of our considerations in~\zita{MR2570113}.
 We are striving for a natural classification of the set of truncated matricial Hamburger moment problems of type ``\(\lleq\)''.
 From \rthm{T1615} we see that these problems have a solution if and only if the sequence of data is \tHnnd{}.
 This leads us to the following relation in the set \(\Hggq{2n}\).

\bnotal{H.N.mpeqR}
 If \(n\in\NO\) and if \(\set{\seqs{2n},\seqt{2n}}\subseteq\Hggq{2n}\), then we write \(\seqs{2n}\mpeqR\seqt{2n}\) if \(\MggqRskg{2n}=\Mggoaakg{q}{\R}{\seqt{2n}}\).
\enota

\breml{H.R.Req}
 Let \(n\in\NO\).
 Then the relation \(\mpeqR\) is an equivalence relation on the set \(\Hggq{2n}\).
\erem

 Let \(n\in\NO\).
 If \(\seqs{2n}\in\Hggq{2n}\), then let \(\mpeqRc{\seqs{2n}}\defeq\setaca{\seqt{2n}\in\Hggq{2n}}{\seqt{2n}\mpeqR\seqs{2n}}\).
 Furthermore, if \(\mathcal{S}\) is a subset of \(\Hggq{2n}\), then let \(\mpeqRc{\mathcal{S}}\defeq\setaca{\mpeqRc{\seqs{2n}}}{\seqs{2n}\in\mathcal{S}}\).

 Looking back to \rthm{T1456} we see that each equivalence class contains a unique representative belonging to \(\Hggeq{2n}\).
 The considerations of \rsec{S.H} provide us now not only detailed insights into the explicit structure of this distinguished representative but even an alternative approach.
 The following notion is the central object of this section.
 
\bdefnl{H.D.Heqseq}
 If $n\in\NO $ and $\seqs{2n}\in\Hggq{2n}$, then the sequence $\Reqseq{s}{2n}$ given by \(\Reqseqmat{s}{2n}\defeq\RLu{n}\), where \(\RLu{n}\) is given in \rnota{H.N.F}, and by \(\Reqseqmat{s}{j}\defeq\su{j}\) for all \(j\in\mn{0}{2n-1}\) is called the \emph{\tReqseq{\seqs{2n}}}.
\edefn

 Now we derive the announced sharpened version of \rthm{T1456}.

 In the following, we will use the notation given in \rdefn{H.D.Heqseq}.

\bpropl{H.L2118}
 Let \(n\in\NO\) and let \(\seqs{2n}\in\Hggq{2n}\).
 Then \(\Reqseq{s}{2n}\in\Hggeq{2n}\) and \(\Reqseq{s}{2n}\mpeqR\seqs{2n}\).
\eprop
\bproof
 In the case \(n=0\), we have \(\Reqseqmat{s}{0}=\RLu{0}=\su{0}\) and \(\Hggq{0}=\Hggeq{0}\).
 Now assume \(n\geq1\).
 According to \rthm{sa4ab5}, we have \(\Reqseqmat{s}{2n}\lleq\su{2n}\) and, in view of \rdefn{H.D.Heqseq}, furthermore, \(\Reqseqmat{s}{j}=\su{j}\) for all \(j\in\mn{0}{2n-1}\).
 From \rrem{H.R1035}, we get then \(\MggqRakg{\Reqseqsy{s}}{2n}\subseteq\MggqRskg{2n}\).
 Conversely, we consider now an arbitrary \(\sigma\in\MggqRskg{2n}\).
 Let the sequence \(\seqa{u}{2n}\) be given by \(u_j\defeq\int_\R x^j\sigma\rk{\dif x}\).
 Then \(\su{2n}-u_{2n}\in\Cggq\) and \(u_j=\su{j}\) for all \(j\in\mn{0}{2n-1}\).
 Furthermore, \(\sigma\in\MggqRag{u}{2n}\), implying \(\seqa{u}{2n}\in\Hggeq{2n}\), by virtue of \rthm{T1708}.
 Consequently, \(u_{2n}\in\exHggesk{2n-1}{\su{2n}}\).
 According to \rthm{sa4ab5}, thus \(u_{2n}\in\matint{\Tripu{n}}{\Reqseqmat{s}{2n}}\).
 In particular, \(\Reqseqmat{s}{2n}-u_{2n}\in\Cggq\).
 Since \(u_j=\su{j}=\Reqseqmat{s}{j}\) for all \(j\in\mn{0}{2n-1}\), then \rrem{H.R1035} yields \(\MggqRakg{u}{2n}\subseteq\MggqRakg{\Reqseqsy{s}}{2n}\).
 Taking additionally into account \rrem{H.R1032}, we can conclude \(\sigma\in\MggqRakg{\Reqseqsy{s}}{2n}\).
 Hence, \(\MggqRskg{2n}\subseteq\MggqRakg{\Reqseqsy{s}}{2n}\).
 Consequently, \(\MggqRakg{\Reqseqsy{s}}{2n}=\MggqRskg{2n}\), implying \(\Reqseq{s}{2n}\mpeqR\seqs{2n}\).
\eproof

\bpropl{H.L3220}
 Let \(n\in\NO\) and let \(\seqs{2n}\in\Hggq{2n}\).
 If \(\seqt{2n}\in\Hggeq{2n}\) satisfies \(\seqt{2n}\mpeqR\seqs{2n}\), then \(\seqt{2n}\) coincides with the \tReqseq{\seqs{2n}}.
\eprop
\bproof
 Let \(\seqt{2n}\in\Hggeq{2n}\) be such that \(\seqt{2n}\mpeqR\seqs{2n}\).
 Observe that \(\seqt{2n}\in\Hggq{2n}\), by virtue of \rrem{H.R.He<gg}.
 In view of \rrem{H.R.Req}, we infer from \rprop{H.L2118} then \(\seqt{2n}\mpeqR\Reqseq{s}{2n}\), \tie{}, \(\MggqRakg{t}{2n}=\MggqRakg{\Reqseqsy{s}}{2n}\).
 According to \rthm{T1708}, we can chose a measure \(\tau\in\MggqRag{t}{2n}\).
 By virtue of \rrem{H.R1032}, then \(\tau\in\MggqRakg{t}{2n}\).
 Thus, \(\tau\in\MggqRakg{\Reqseqsy{s}}{2n}\).
 Consequently, we have \(t_{2n}=\int_\R x^{2n}\tau\rk{\dif x}\lleq\Reqseqmat{s}{2n}\) and \(t_j=\int_\R x^j\tau\rk{\dif x}=\Reqseqmat{s}{j}\) for all \(j\in\mn{0}{2n-1}\).
 Since \(\Reqseq{s}{2n}\) belongs to \(\Hggeq{2n}\), according to \rprop{H.L2118}, we can conclude in a similar way \(\Reqseqmat{s}{2n}\lleq t_{2n}\).
 Hence, \(t_{2n}=\Reqseqmat{s}{2n}\) follows.
\eproof

 Now we state the main result of this section, which sharpens \rthm{T1456}.

\bthml{H.P0828}
 Let \(n\in\NO\) and let \(\seqs{2n}\in\Hggq{2n}\).
 Then \(\mpeqRc{\seqs{2n}}\cap\Hggeq{2n}=\set{\Reqseq{s}{2n}}\).
\ethm
\bproof
 Combine \rpropss{H.L2118}{H.L3220}.
\eproof
 
 Our next aim can be described as follows.
 Let \(n\in\N\) and let \(\seqs{2n}\in\Hggq{2n}\).
 Then an appropriate application of \rthm{A.T1514} leads us to the determination of all sequences \(\seqa{r}{2n}\) which are contained in the equivalence class \(\mpeqRc{\seqs{2n}}\).
 
\bpropl{H.P0830}
 Let \(n\in\N\) and let \(\seqs{2n}\in\Hggq{2n}\).
 Then \(\mpeqRc{\seqs{2n}}\) coincides with the set of all sequences \(\seqa{r}{2n}\) of complex \tqqa{matrices} fulfilling \(\ran{r_{2n}-\RLu{n}}\cap\ran{\Lu{n-1}}=\set{\Ouu{q}{1}}\), \(r_{2n}-\RLu{n}\in\Cggq\), and \(r_j=\su{j}\) for all \(j\in\mn{0}{2n-1}\).
\eprop
\bproof
 Let \(\seqt{2n}\) be the \tReqseq{\seqs{2n}}.
 By virtue of \rprop{H.L2118}, we have \(\seqt{2n}\in\Hggeq{2n}\) and \(\seqt{2n}\mpeqR\seqs{2n}\).
 According to \rdefn{H.D.Heqseq}, furthermore \(t_{2n}=\RLu{n}\) and \(t_{j}=\su{j}\) for all \(j\in\mn{0}{2n-1}\).
 In particular, \(\Luo{n-1}{t}=\Lu{n-1}\) and \(\Tripuo{n}{t}=\Tripu{n}\).
 
 Consider now an arbitrary \(\seqa{r}{2n}\in\mpeqRc{\seqs{2n}}\), \tie{}, \(\seqa{r}{2n}\in\Hggq{2n}\) with \(\seqa{r}{2n}\mpeqR\seqs{2n}\).
 In particular, \(\seqt{2n}\mpeqR\seqa{r}{2n}\), by virtue of \rrem{H.R.Req}.
 From \rprop{H.L3220} we can conclude then that \(\seqt{2n}\) coincides with the  \tReqseq{\seqa{r}{2n}}.
 In view of \rdefn{H.D.Heqseq}, consequently \(t_{2n}=\RLuo{n}{r}\) and \(t_{j}=r_j\) for all \(j\in\mn{0}{2n-1}\).
 Hence, \(\RLu{n}=\RLuo{n}{r}\) and \(\su{j}=r_j\) for all \(j\in\mn{0}{2n-1}\).
 In particular, \(\Lu{n-1}=\Luo{n-1}{r}\) and \(\Tripu{n}=\Tripuo{n}{r}\).
 From \rrem{H.R.L>} we know that \(\Luo{n}{r}\in\Cggq\).
 Setting  \(A\defeq\Luo{n}{r}\) and \(\cV\defeq\ran{\Luo{n-1}{r}}\), we obtain, in view of \rnota{H.N.F} and \rnotap{H.N.yz}{H.N.yz.c}, then
\[
 r_{2n}-\RLu{n}
 =r_{2n}-\RLuo{n}{r}
 =r_{2n}-\Tripuo{n}{r}-\SO{A}{\cV}
 =\Luo{n}{r}-\SO{A}{\cV}
 =A-\SO{A}{\cV}.
\]
 Because of \(\ran{\Lu{n-1}}=\cV\) and \rthmp{A.T1514}{A.T1514.a}, thus \(\ran{r_{2n}-\RLu{n}}\cap\ran{\Lu{n-1}}=\set{\Ouu{q}{1}}\).
 Furthermore, \rpropp{sa1ab3}{sa1ab3.a} yields \(r_{2n}-\RLu{n}\in\Cggq\).
 
 Conversely, we consider now an arbitrary sequence \(\seqa{r}{2n}\) of complex \tqqa{matrices} fulfilling \(\ran{r_{2n}-\RLu{n}}\cap\ran{\Lu{n-1}}=\set{\Ouu{q}{1}}\), \(r_{2n}-\RLu{n}\in\Cggq\), and \(r_j=\su{j}\) for all \(j\in\mn{0}{2n-1}\).
 Then \(\Luo{n-1}{r}=\Lu{n-1}\) and \(\Tripuo{n}{r}=\Tripu{n}\).
 Because of \(\seqs{2n}\in\Hggq{2n}\), we have \(\seqa{r}{2n-1}\in\Hggeq{2n-1}\).
 From \rlem{H.L0901} we infer that \(\RLu{n}^\ad=\RLu{n}\) with \(\Tripu{n}\lleq\RLu{n}\).
 Consequently, \(r_{2n}\in\CHq\) and \(\Tripu{n}\lleq\RLu{n}\lleq r_{2n}\).
 Hence, \(\Tripuo{n}{r}\lleq r_{2n}\), \tie{}, \(\Luo{n}{r}\in\Cggq\).
 By virtue of \rprop{H.P0840}, we obtain then \(\seqa{r}{2n}\in\Hggq{2n}\).
 Denote by $\Reqseq{r}{2n}$ the \tReqseq{\seqa{r}{2n}} and let \(A\defeq\Luo{n}{r}\) and \(\cV\defeq\ran{\Luo{n-1}{r}}\).
 By \rdefn{H.D.Heqseq} and \rnota{H.N.F}, then \(\Reqseqmat{r}{2n}=\Tripuo{n}{r}+\SO{A}{\cV}\) and \(\Reqseqmat{r}{j}=r_{j}\) for all \(j\in\mn{0}{2n-1}\).
 Consequently, \(\Reqseqmat{r}{j}=r_{j}=\su{j}=t_j\) for all \(j\in\mn{0}{2n-1}\).
 Setting \(X\defeq t_{2n}-\Tripuo{n}{r}\) and  \(Y\defeq r_{2n}-t_{2n}\), we have, by virtue of \rnotap{H.N.yz}{H.N.yz.c}, then \(X+Y=r_{2n}-\Tripuo{n}{r}=\Luo{n}{r}=A\) and \(X=t_{2n}-\Tripu{n}=t_{2n}-\Tripuo{n}{t}=\Luo{n}{t}\) and, furthermore, \(Y=r_{2n}-\RLu{n}\).
 By assumption, then \(\ran{Y}\cap\ran{\Lu{n-1}}=\set{\Ouu{q}{1}}\) and \(Y\in\Cggq\).
 From \rrem{H.R.He<gg} we infer \(\seqt{2n}\in\Hggq{2n}\).
 In particular, \(\seqt{2n-1}\in\Hggeq{2n-1}\) and \(t_{2n}\in\exHggs{2n-1}\).
 \rprop{sa2ab5} yields then \(\Luo{n}{t}\in\Cggq\) and \(\ran{\Luo{n}{t}}\subseteq\ran{\Luo{n-1}{t}}\).
 Hence, \(X\in\Cggq\).
 Taking into account \(\Luo{n-1}{t}=\Lu{n-1}=\Luo{n-1}{r}\), we see furthermore \(\ran{X}\subseteq\cV\) and \(\ran{Y}\cap\cV=\set{\Ouu{q}{1}}\).
 From \rthmp{A.T1514}{A.T1514.b} we get then \(X=\SO{A}{\cV}\).
 Hence, \(\Reqseqmat{r}{2n}=\Tripuo{n}{r}+\SO{A}{\cV}=t_{2n}\).
 Thus, the sequences $\Reqseq{r}{2n}$  and \(\seqt{2n}\) coincide.
 Using \rprop{H.L2118} and \rrem{H.R.Req}, we get then \(\seqa{r}{2n}\mpeqR\Reqseq{r}{2n}\mpeqR\seqt{2n}\mpeqR\seqs{2n}\).
 Consequently, \(\seqa{r}{2n}\in\mpeqRc{\seqs{2n}}\).
\eproof

\section{On truncated matricial $\rhl$-Stieltjes moment problems}\label{S.KMP}
 In our following considerations, let \(\ug\) be a real number.
 In order to state a necessary and sufficient condition for the solvability of each of the moment problems \mprob{\rhl}{m}{\lleq} and \mprob{\rhl}{m}{=}, we have to recall the notion of two types of sequences of matrices.

 Let \(\kappa\in\Ninf\) and let \(\seqska\) is a sequence of complex \tpqa{matrices}.
 For each \(n\in\NO\) with \(2n+1\leq\kappa\), let the block Hankel  matrix \(\Ku{n}\) be given by \(\Ku{n}\defeq\matauuuo{s_{j+k+1}}{j}{k}{0}{n}\).
 Furthermore, let the sequence \(\seqsa{\kappa-1}\) be given by \(\sa{j}\defeq-\ug\su{j}+\su{j+1}\).
 For each matrix \(X_k=X_k^{\rk{s}}\) built from the sequence \(\seqska\), denote (if possible) by \(X_{\ug,k}\defeq X_k^{\rk{\saus}}\) the corresponding matrix built from the sequence \(\seqsa{\kappa-1}\) instead of \(\seqska\).
 In particular, we have then \(\Hau{n}=-\ug\Hu{n}+\Ku{n}\) for all \(n\in\NO\) with \(2n+1\leq\kappa\).
 In the classical case \(\ug=0\), we see that \(\sa{j}=\su{j+1}\) for all \(j\in\mn{0}{\kappa-1}\).

 Let \(\Kggq{0}\defeq\Hggq{0}\).
 For each \(n\in\N\), denote by \(\Kggq{2n}\) the set of all sequences \(\seqs{2n}\) of complex \tqqa{matrices} for which the \tbHms{} \(\Hu{n}\) and \(\Hau{n-1}\) are both \tnnH{}.
 For each \(n\in\NO\), denote by \(\Kggq{2n+1}\) the set of all sequences \(\seqs{2n+1}\) of complex \tqqa{matrices} for which the \tbHms{} \(\Hu{n}\) and \(\Hau{n}\) are both \tnnH{}.
 Furthermore, denote by \(\Kggqinf\) the set of all sequences \(\seqsinf\) of complex \tqqa{matrices} satisfying \(\seqs{m}\in\Kggq{m}\) for all \(m\in\NO\).
 The sequences belonging to \(\Kggq{0}\), \(\Kggq{2n}\), \(\Kggq{2n+1}\), or \(\Kggqinf\) are said to be \emph{\tKnnd}.

 Now we can characterize the situations that the mentioned problems have a solution:

\bthmnl{\zitaa{MR2735313}{\cthm{1.4}{909}}}{dm2.T11}
 Let \(m\in\NO\) and let \(\seqs{m}\) be a sequence of complex \tqqa{matrices}.
 Then \(\MggqKskg{m}\neq\emptyset\) if and only if \(\seqs{m}\in\Kggq{m}\).
\ethm

 Now we characterize the solvability of Problem~\mprob{\rhl}{m}{=}.

\bthmnl{\zitaa{MR2735313}{\cthm{1.3}{909}}}{K.T1435}
 Let \(m\in\NO\) and let \(\seqs{m}\) be a sequence of complex \tqqa{matrices}.
 Then \(\MggqKsg{m}\neq\emptyset\) if and only if \(\seqs{m}\in\Kggeq{m}\).
\ethm

 The following result is the starting point of our subsequent considerations:
 
\bthmnl{\zitaa{MR2735313}{\cthm{5.2}{935}}}{K.T2547}
 Let \(m\in\NO\) and let \(\seqs{m}\in\Kggq{m}\).
 Then there exists a unique sequence \(\Keqseq{s}{m}\in\Kggeq{m}\) such that
\beql{K.T2547.A}
 \Mggoaakg{q}{\rhl}{\Keqseq{s}{m}}
 =\MggqKskg{m}.
\eeq
\ethm

 \rthm{K.T2547} was very essential for the considerations in~\zita{MR2735313}.
 
 The main goal of the rest of this paper is to derive this result by use of an appropriate application of the machinery developed in \rsec{S.ando}.
 This will lead us to an explicit formula for the desired sequence \(\Keqseq{s}{m}\).
 
 Following~\zita{MR2735313} we sketch now some essential features of the history of \rthm{K.T2547}.
 In the case \(\ug=0\) the existence of a sequence \(\Keqseq{s}{m}\in\Kggeq{m}\) satisfying \eqref{K.T2547.A} was already formulated by V.~A.~Bolotnikov~\zitaa{MR1362524}{\cthm{1.5}{443}, \clem{1.6}{443}}.
 This result is true.
 However, it was shown in~\zitaa{MR2735313}{\cexam{5.1}{934}} that the concrete sequence \(\Keqseq{s}{m}\) constructed in~\zitaa{MR1362524}{\clemss{2.7}{450}{6.3}{466}} does not produce a moment problem equivalent to \mprob{[0,\infp)}{m}{\lleq}.

\section{On a restricted extension problem for a finite \hKnnde{} sequence}\label{S.K}
 This section is written against to the background of applying the results of \rsec{S.ando} to the truncated matricial \(\rhl\)\nobreakdash-Stieltjes moment problems formulated in the introduction for \(\Omega=\rhl\).
 In the heart of our strategy lies the treatment of a special restricted extension problem for matrices.
 The complete answer to this problem is contained in \rthm{K.P1416} which is the central result of this section.

 Using \rrem{H.R0857} we can conclude:
 
\breml{K.R1447}
 Let \(\kappa\in\NOinf\) and let \(\seqska\in\Kggqka\).
 Then \(\su{j}\in\CHq\) for all \(j\in\mn{0}{\kappa}\) and \(\su{2k}\in\Cggq\) for all \(k\in\NO\) with \(2k\leq\kappa\).
 If \(\kappa\geq1\), furthermore \(\sa{j}\in\CHq\) for all \(j\in\mn{0}{\kappa-1}\) and \(\sa{2k}\in\Cggq\) for all \(k\in\NO\) with \(2k\leq\kappa-1\).
\erem

\bdefnnl{\tcf{}~\zitaa{MR3014201}{\cdefn{4.2}{223}}}{K.D.kpf}
 If \(\seqska\) is a sequence of complex \tpqa{matrices}, then the sequence \(\seq{\kpu{j}}{j}{0}{\kappa}\) given by
 \(
  \kpu{2k}
  \defeq\su{2k}-\Tripu{k}
 \)
 for all \(k\in\NO\) with \(2k\leq\kappa\) and by
 \(
  \kpu{2k+1}
  \defeq\sa{2k}-\Tripa{k}
 \)
 for all \(k\in\NO\) with \(2k+1\leq\kappa\) is called the \emph{\tkpfa{\(\seqska \)}}.
\edefn
 
\breml{K.R.Q>}
 If \(\kappa\in\NOinf\) and \(\seqska\in\Kggqka\), then one can easily see from~\zitaa{MR3014201}{\cthm{4.12(b)}{}} that \(\seqs{m}\in\Kggq{m}\) and \(\kpu{m}\in\Cggq\) for all \(m\in\mn{0}{\kappa}\).
\erem

\bnotal{K.N2530}
 If \(m\in\NO\) and \(\seqs{m}\) is a sequence of complex \tqqa{matrices}, then denote by \(\exKggs{m}\) the set of all complex \tqqa{matrices} \(\su{m+1}\) such that \(\seqs{m+1}\in\Kggq{m+1}\).
\enota

 For each \(m\in\NO\), denote by \(\Kggeq{m}\) the set of all sequences \(\seqs{m}\) of complex \tqqa{matrices} for which there exists a complex \tqqa{matrix} \(\su{m+1}\) such that the sequence \(\seqs{m+1}\) belongs to \(\Kggq{m+1}\).
 Furthermore, let \(\Kggeqinf\defeq\Kggqinf\).
 The sequences belonging to \(\Kggeq{m}\) or \(\Kggeqinf\) are said to be \emph{\tKnnde}.
 Obviously, if \(m\in\NO\) and if \(\seqs{m}\) is a sequence of complex \tqqa{matrices}, then \(\exKggs{m}\neq\emptyset\) if and only if \(\seqs{m}\in\Kggeq{m}\)
 
 The following result shows why the class \(\Kggeqka\) is important.

\bthmnl{\zitaa{MR3014201}{\cthm{1.6}{214}}}{dm2.T10}
 Let \(\kappa\in\NOinf\) and let \(\seqska\) be a sequence of complex \tqqa{matrices}.
 Then \(\MggqKsg{\kappa}\neq\emptyset\) if and only if \(\seqska \in\Kggeqka\).
\ethm

\breml{K.R.Ke<gg}
 From \rrem{K.R.Q>} one can easily see that \(\Kggeqka\subseteq\Kggqka\) for all \(\kappa\in\NOinf\).
\erem

\bpropnl{\tcf{}~\zitaa{MR3014201}{\cthm{4.12(b),~(c)}{}}}{K.P1430}
 Let \(m\in\NO\), let \(\seqs{m}\in\Kggeq{m}\), and let \(\su{m+1}\in\Cqq\).
 Then  \(\su{m+1}\in\exKggs{m}\) if and only if \(\kpu{m+1}\in\Cggq\).
\eprop

\bnotal{K.N1439}
 If \(\kappa\in\NOinf\) and \(\seqska\) is a sequence of complex \tpqa{matrices}, then let \(\umg{2k-1}\defeq\Tripu{k}\) for all \(k\in\NO\) with \(2k-1\leq\kappa\) and let \(\umg{2k}\defeq\ug\su{2k}+\Tripa{k}\) for all \(k\in\NO\) with \(2k\leq\kappa\).
\enota

\breml{K.R.Q=s-a}
 If \(\kappa\in\NOinf\) and \(\seqska\) is a sequence of complex \tpqa{matrices}, then \rdefn{K.D.kpf} shows that \(\kpu{j}=\su{j}-\umg{j-1}\) for all \(j\in\mn{0}{\kappa}\).
\erem

\breml{K.R1509}
 If \(\kappa\in\NOinf\) and \(\seqska\in\Kggeqka\), then \rremss{K.R.Q>}{K.R1447} show that \(\umg{m}^\ad=\umg{m}\) for all \(m\in\mn{-1}{\kappa}\).
\erem

\bnotal{K.N1435}
 Let \(m\in\NO\), let \(\seqs{m}\) be a sequence of complex \tqqa{matrices}, and let \(Y\in\CHq\).
 Then denote by \(\exKggsk{m}{Y}\) the set of all \(X\in\exKggs{m}\) satisfying \(Y-X\in\Cggq\).
\enota

\bleml{K.L1437}
 If \(\kappa\in\NOinf\) and \(\seqska\in\Kggqka\), then \(\umg{m-1}\lleq\su{m}\) for all \(m\in\mn{0}{\kappa}\).
\elem
\bproof
 Obviously, \(\umg{-1}=\Oqq\lleq\su{0}\) by virtue of \rrem{K.R1447}.
 Now assume \(\kappa\geq1\) and we consider an arbitrary \(m\in\mn{1}{\kappa}\).
 \rrem{K.R1447} yields \(\su{m}\in\CHq\).
 \rrem{K.R.Q>} shows that \(\kpu{m}\in\Cggq\).
 In view of \rrem{K.R.Q=s-a}, then \(\umg{m-1}\in\CHq\) and \(\umg{m-1}\lleq\su{m}\).
\eproof

\bthml{K.P1517}
 If \(m\in\N\) and \(\seqs{m}\in\Kggq{m}\), then \(\umg{m-1}\lleq\su{m}\) and \(\exKggsk{m-1}{\su{m}}=\matint{\umg{m-1}}{\su{m}}\).
\ethm
\bproof
 By virtue of \rlem{K.L1437}, we have \(\set{\umg{m-1},\su{m}}\subseteq\CHq\) and \(\umg{m-1}\lleq\su{m}\).
 Let the sequence \(\seq{t_j}{j}{0}{m-1}\) be given by \(t_j\defeq\su{j}\) for all \(j\in\mn{0}{m-1}\).
 Obviously, then \(\umgo{m-1}{t}=\umg{m-1}\).
 Furthermore, \(\set{\seqs{m-1},\seqt{m-1}}\subseteq\Kggeq{m-1}\).
 
 First consider now an arbitrary \(t_{m}\in\exKggsk{m-1}{\su{m}}\).
 Then \(\seqt{m}\in\Kggq{m}\) and \(\su{m}-t_{m}\in\Cggq\).
 In particular, \(t_{m}\in\CHq\).
 In view of \rrem{K.R.Q=s-a}, \(t_{m}-\umg{m-1}=t_{m}-\umgo{m-1}{t}=\kpuo{m}{t}\), and \rrem{K.R.Q>}, we get then \(\umg{m-1}\lleq t_{m}\lleq\su{m}\).
 
 Conversely, let \(t_{m}\in\matint{\umg{m-1}}{\su{m}}\).
 Then \(\su{m}-t_{m}\in\Cggq\) and, in view of \rrem{K.R.Q=s-a}, furthermore \(\kpuo{m}{t}=t_{m}-\umgo{m-1}{t}=t_{m}-\umg{m-1}\).
 In particular, \(\kpuo{m}{t}\) is \tnnH{}.
 According to \rprop{K.P1430}, thus \(\seqt{m}\in\Kggq{m}\).
 Consequently, \(t_{m}\in\exKggsk{m-1}{\su{m}}\).
\eproof

\bcorl{H.C1552}
 If \(m\in\NO\) and \(\seqs{m}\in\Kggq{m}\), then \(\setaca{\mpm{\sigma}{m}}{\sigma\in\MggqKskg{m}}\subseteq\matint{\umg{m-1}}{\su{m}}\).
\ecor
\bproof
 In view of \(\umg{-1}=\Tripu{0}=\Oqq\), the case \(m=0\) is obvious.
 If \(m\geq1\), combine \rthm{dm2.T10}, \rrem{K.R.Ke<gg}, and \rthm{K.P1517}.
\eproof

\bnotal{K.N1010}
 If \(m\in\NO\) and \(\seqs{m}\) is a sequence of complex \tqqa{matrices}, then denote by \(\exKgges{m}\) the set of all complex \tqqa{matrices} \(\su{m+1}\) such that \(\seqs{m+1}\in\Kggeq{m+1}\).
\enota

 The following result is essential for the realization of our concept of a new approach to \rthm{K.T2547}.

\bthmnl{\tcf{}~\zitaa{MR3014201}{\cthm{4.12(c)}{}}}{K.P1016}
 Let \(m\in\NO\), let \(\seqs{m}\in\Kggeq{m}\), and let \(\su{m+1}\in\Cqq\).
 Then  \(\su{m+1}\in\exKgges{m}\) if and only if \(\kpu{m+1}\in\Cggq\) and 
\(\ran{\kpu{m+1}}\subseteq\ran{\kpu{m}}\).
\ethm

\bnotal{K.N1443}
 Let \(m\in\NO\), let \(\seqs{m}\) be a sequence of complex \tqqa{matrices}, and let \(Y\in\CHq\).
 Then denote by \(\exKggesk{m}{Y}\) the set of all \(X\in\exKgges{m}\) satisfying \(Y-X\in\Cggq\).
\enota

 Observe that the following construction is well defined due to \rrem{K.R.Q>}:
 
\bnotal{K.N.grG}
 If \(\kappa\in\NOinf\) and \(\seqska\in\Kggqka\), then let \(\Rku{0}\defeq\su{0}\) and let \(\Rku{m}\defeq\umg{m-1}+\SO{\kpu{m}}{\ran{\kpu{m-1}}}\) for all \(m\in\mn{1}{\kappa}\).
\enota

\bleml{K.L1840}
 Let \(\kappa\in\NOinf\) and let \(\seqska\in\Kggqka\).
 Then \(\Rku{m}^\ad=\Rku{m}\) and \(\umg{m-1}\lleq\Rku{m}\lleq\su{m}\) for all \(m\in\mn{0}{\kappa}\).
\elem
\bproof
 Let \(m\in\mn{0}{\kappa}\).
 From \rlem{K.L1437} we see that the matrices \(\umg{m-1}\) and \(\su{m}\) are \tH{} and that \(\umg{m-1}\lleq\su{m}\).
 \rrem{N23} shows that \(\Rku{m}^\ad=\Rku{m}\).
 By virtue of \(\Rku{0}=\su{0}\), we get in particular \(\umg{-1}\lleq\Rku{0}\lleq\su{0}\).
 Now assume \(\kappa\geq1\) and consider an arbitrary \(m\in\mn{1}{\kappa}\).
 According to \rrem{K.R.Q>}, we have \(\kpu{m}\in\Cggq\).
 From \rrem{N23} and \rpropp{sa1ab3}{sa1ab3.a} we obtain \(\Oqq\lleq\SO{\kpu{m}}{\ran{\kpu{m-1}}}\lleq\kpu{m}\), implying, by virtue of \rnota{K.N.grG} and \rrem{K.R.Q=s-a}, then \(\umg{m-1}\lleq\Rku{m}\lleq\su{m}\).
\eproof

\bthml{K.P1416}
 Let \(m\in\N\) and let \(\seqs{m}\in\Kggq{m}\).
 Then:
\benui
 \il{K.P1416.a} \(\umg{m-1}\lleq\Rku{m}\lleq\su{m}\) and \(\exKggesk{m-1}{\su{m}}=\matint{\umg{m-1}}{\Rku{m}}\).
 \il{K.P1416.b} \(\Rku{m}=\su{m}\) if and only if \(\seqs{m}\in\Kggeq{m}\).
\eenui
\ethm
\bproof
 \eqref{K.P1416.a} From \rlem{K.L1840}, we see that the matrices \(\umg{m-1}\), \(\Rku{m}\), and \(\su{m}\) are \tH{} and that \(\umg{m-1}\lleq\Rku{m}\lleq\su{m}\).
 In view of \rrem{K.R.Q>}, we get \(\kpu{m}\in\Cggq\).
 Setting \(A\defeq\kpu{m}\) and \(\cV\defeq\ran{\kpu{m-1}}\), we have \(\Rku{m}=\umg{m-1}+\SO{A}{\cV}\) and \(\su{m}=\umg{m-1}+A\), by virtue of \rrem{K.R.Q=s-a}.
 Let the sequence \(\seq{t_j}{j}{0}{m-1}\) be given by \(t_j\defeq\su{j}\).
 Then \(\umgo{m-1}{t}=\umg{m-1}\) and \(\kpuo{m-1}{t}=\kpu{m-1}\).
 In particular, \(\cV=\ran{\kpuo{m-1}{t}}\).
 Furthermore, \(\set{\seqs{m-1},\seq{t_j}{j}{0}{m-1}}\subseteq\Kggeq{m-1}\).
 
 First consider an arbitrary \(t_{m}\in\exKggesk{m-1}{\su{m}}\).
 Then \(\seqt{m}\in\Kggeq{m}\), \(t_{m}\in\CHq\), and \(\su{m}-t_{m}\in\Cggq\).
 Let \(X\defeq\kpuo{m}{t}\).
 According to \rthm{K.P1016}, then  \(X\in\Cggq\) and \(\ran{X}\subseteq\cV\).
 In view of \rrem{K.R.Q=s-a}, we have \(X=t_{m}-\umgo{m-1}{t}=t_{m}-\umg{m-1}\) and, hence, \(X\leq\su{m}-\umg{m-1}=A\) follows.
 Consequently, \(X\in\LcR{A}{\cV}\).
 The application of \rthm{sa3ab3} yields then \(\SO{A}{\cV}\lgeq X\).
 In view of \rnota{K.N.grG}, thus \(\umg{m-1}\lleq t_{m}\lleq\Rku{m}\).
 
 Conversely, let \(t_{m}\in\matint{\umg{m-1}}{\Rku{m}}\), \tie{}, \(t_{m}\in\CHq\) with \(\umg{m-1}\lleq t_{m}\lleq\Rku{m}\).
 Then \(X\defeq t_{m}-\umg{m-1}\) is \tH{} and fulfills \(\Oqq\lleq X\lleq\SO{A}{\cV}\).
 By virtue of \rrem{A.R1004}, then \(\ran{X}\subseteq\ran{\SO{A}{\cV}}\) and, according to \rprop{sa1ab3}, furthermore \(\ran{\SO{A}{\cV}}=\ran{A}\cap\cV\) and \(\SO{A}{\cV}\lleq A\).
 Consequently, \(\ran{X}\subseteq\cV\) and \(\Oqq\lleq X\lleq A\).
 In view of \rrem{K.R.Q=s-a}, furthermore \(\kpuo{m}{t}=t_{m}-\umgo{m-1}{t}=t_{m}-\umg{m-1}=X\).
 Therefore, \(\kpuo{m}{t}\in\Cggq\) and \(\ran{\kpuo{m}{t}}\subseteq\ran{\kpuo{m-1}{t}}\).
 According to \rthm{K.P1016}, thus \(t_m\in\exKgge{t}{m-1}\), \tie{}, \(\seqt{m}\in\Kggeq{m}\).
 Moreover, since \rrem{K.R.Q=s-a} shows that \(\su{m}-t_{m}=A-X\) we see that \(\su{m}-t_m\in\Cggq\).
 Consequently, \(t_{m}\in\exKggesk{m-1}{\su{m}}\).
 
 \eqref{K.P1416.b} In view of \(\seqs{m}\in\Kggq{m}\), we have \(\seqs{m-1}\in\Kggeq{m-1}\) and \(\su{m}\in\exKggs{m-1}\).
 \rprop{K.P1430} yields then \(\kpu{m}\in\Cggq\).
 By virtue of \rnota{K.N.grG} and \rrem{K.R.Q=s-a} we have \(\Rku{m}=\su{m}\) if and only if \(\kpu{m}=\SO{\kpu{m}}{\ran{\kpu{m-1}}}\).
 According to \rpropp{sa1ab3}{sa1ab3.g}, the latter is equivalent to \(\ran{\kpu{m}}\subseteq\ran{\kpu{m-1}}\).
 In view of \(\seqs{m-1}\in\Kggeq{m-1}\) and \(\kpu{m}\in\Cggq\), \rthm{K.P1016} shows the equivalence of \(\ran{\kpu{m}}\subseteq\ran{\kpu{m-1}}\) and \(\su{m}\in\exKgges{m-1}\).
 Hence,~\eqref{K.P1416.b} follows.
\eproof

 \rthm{K.P1416} leads us now quickly in an alternative way to one of the main results of~\zita{MR2735313}.

\bthmnl{\tcf{}~\zitaa{MR2735313}{\cthm{5.4}{939}}}{K.C1544}
 If \(m\in\NO\) and \(\seqs{m}\in\Kggq{m}\), then \(\setaca{\mpm{\sigma}{m}}{\sigma\in\MggqKskg{m}}=\matint{\umg{m-1}}{\Rku{m}}\).
\ethm
\bproof
 In view of \(\umg{-1}=\Tripu{0}=\Oqq\) and \(\Rku{0}=\su{0}\), the case \(m=0\) is obvious.
 If \(m\geq1\), combine \rthmss{dm2.T10}{K.P1416}.
\eproof

\section{On equivalence classes of truncated matricial moment problems of type \mprob{\rhl}{m}{\lleq}}\label{S.Keq}
 This section contains an aspect of our considerations in~\zita{MR2735313}.
 We are striving for a natural classification of the set of truncated matricial \(\rhl\)\nobreakdash-Stieltjes moment problems of type ``\(\lleq\)''.
 From \rthm{dm2.T11} we see that these problems have a solution if and only if the sequence of data is \tKnnd{}.
 This leads us to the following relation in the set \(\Kggq{m}\).

 If \(m\in\NO\) and \(\seqs{m},\seqt{m}\) are two sequences belonging to \(\Kggq{m}\), then we write \(\seqs{m}\mpeqK\seqt{m}\) if \(\MggqKskg{m}=\MggqKakg{t}{m}\) holds true.

\breml{K.R.Keq}
 It is readily checked that \(\mpeqK\) is an equivalence relation on the set \(\Kggq{m}\).
\erem

\bnotal{K.N.mpeqKc}
 Let \(m\in\NO\).
 If \(\seqs{m}\in\Kggq{m}\), then let \(\mpeqKc{\seqs{m}}\defeq\setaca{\seqt{m}\in\Kggq{m}}{\seqt{m}\mpeqK\seqs{m}}\).
 Furthermore, if \(\mathcal{S}\) is a subset of \(\Kggq{m}\), then let \(\mpeqKc{\mathcal{S}}\defeq\setaca{\mpeqKc{\seqs{m}}}{\seqs{m}\in\mathcal{S}}\).
\enota

 Regarding \zitaa{MR2735313}{\cthm{5.2}{935}} we see that each equivalence class contains a unique representative belonging to \(\Kggeq{m}\).
 The considerations of \rsec{S.K} provide us now not only detailed insights into the explicit structure of this distinguished representative but even an alternative approach.
 The following notion is the central object of this section.
 
\bdefnl{K.D.Keqseq}
 If $m\in\NO $ and $\seqs{m}\in\Kggq{m}$, then the sequence $\Keqseq{s}{m}$ given by \(\Keqseqmat{s}{m}\defeq\Rku{m}\), where \(\Rku{m}\) is given in \rnota{K.N.grG}, and by \(\Keqseqmat{s}{j}\defeq\su{j}\) for all \(j\in\mn{0}{m-1}\) is called the \emph{\tKeqseq{\seqs{m}}}.
\edefn

 Now we derive a sharpened version of \zitaa{MR2735313}{\cthm{5.2}{935}}.

\bpropl{K.L2118}
 Let \(m\in\NO\) and let \(\seqs{m}\in\Kggq{m}\).
 Then the \tKnnde{} sequence \(\Keqseq{s}{m}\) equivalent to \(\seqs{m}\) belongs to \(\Kggeq{m}\) and \(\Keqseq{s}{m}\mpeqK\seqs{m}\).
\eprop
\bproof
 Let \(m\in\NO\) and let \(\seqs{m}\in\Kggq{m}\).
 First assume \(m=0\).
 We have \(\Keqseqmat{s}{0}=\Rku{0}=\su{0}\).
 In view of \(\Kggq{0}=\Kggeq{0}\), then the assertions follow.
 Now assume \(m\geq1\).
 According to \rthm{K.P1416}, we have \(\Keqseqmat{s}{m}\lleq\su{m}\) and, in view of \rdefn{K.D.Keqseq}, furthermore \(\Keqseqmat{s}{j}=\su{j}\) for all \(j\in\mn{0}{m-1}\).
 From \rrem{H.R1035}, we get then \(\MggqKakg{\Keqseqsy{s}}{m}\subseteq\MggqKskg{m}\).
 Conversely, let \(\sigma\in\MggqKskg{m}\).
 Let the sequence \(\seqa{u}{m}\) be given by \(u_j\defeq\int_\rhl x^j\sigma\rk{\dif x}\).
 Then \(\su{m}-u_{m}\in\Cggq\) and \(u_j=\su{j}\) for all \(j\in\mn{0}{m-1}\).
 Furthermore, \(\sigma\in\MggqKag{u}{m}\), implying \(\seqa{u}{m}\in\Kggeq{m}\), by virtue of \rthm{dm2.T10}.
 Consequently, \(u_{m}\in\exKggesk{m-1}{\su{m}}\).
 According to \rthm{K.P1416}, thus \(u_{m}\in\matint{\umg{m-1}}{\Keqseqmat{s}{m}}\).
 In particular, \(\Keqseqmat{s}{m}-u_{m}\in\Cggq\).
 Since \(u_j=\su{j}=\Keqseqmat{s}{j}\) for all \(j\in\mn{0}{m-1}\), then \rrem{H.R1035} yields \(\MggqKakg{u}{m}\subseteq\MggqKakg{\Keqseqsy{s}}{m}\).
 Taking additionally into account \rrem{H.R1032}, we can conclude \(\sigma\in\MggqKakg{\Keqseqsy{s}}{m}\).
 Hence, we have shown \(\MggqKskg{m}\subseteq\MggqKakg{\Keqseqsy{s}}{m}\).
 Consequently, \(\MggqKakg{\Keqseqsy{s}}{m}=\MggqKskg{m}\), implying \(\Keqseq{s}{m}\mpeqK\seqs{m}\).
\eproof

\bpropl{K.L3220}
 Let \(m\in\NO\) and let \(\seqs{m}\in\Kggq{m}\).
 If \(\seqt{m}\in\Kggeq{m}\) satisfies \(\seqt{m}\mpeqK\seqs{m}\), then \(\seqt{m}\) coincides with the \tKeqseq{\seqs{m}}.
\eprop
\bproof
 Let \(\seqt{m}\in\Kggeq{m}\) be such that \(\seqt{m}\mpeqK\seqs{m}\).
 Observe that \(\seqt{m}\in\Kggq{m}\), by virtue of \rrem{K.R.Ke<gg}.
 In view of \rrem{K.R.Keq}, we infer from \rprop{K.L2118} then \(\seqt{m}\mpeqK\Keqseq{s}{m}\), \tie{}\ \(\MggqKakg{t}{m}=\MggqKakg{\Keqseqsy{s}}{m}\).
 According to \rthm{dm2.T10}, we can chose a measure \(\tau\in\MggqKag{t}{m}\).
 By virtue of \rrem{H.R1032}, then \(\tau\in\MggqKakg{t}{m}\).
 Thus, \(\tau\in\MggqKakg{\Keqseqsy{s}}{m}\).
 Consequently, we have \(t_{m}=\int_\rhl x^{m}\tau\rk{\dif x}\lleq\Keqseqmat{s}{m}\) and \(t_j=\int_\rhl x^j\tau\rk{\dif x}=\Keqseqmat{s}{j}\) for all \(j\in\mn{0}{m-1}\).
 Since \(\Keqseq{s}{m}\) belongs to \(\Kggeq{m}\), according to \rprop{K.L2118}, we can conclude in a similar way \(\Keqseqmat{s}{m}\lleq t_{m}\).
 Hence, \(t_{m}=\Keqseqmat{s}{m}\) follows.
\eproof

\bthml{K.P0828}
 If \(m\in\NO\) and \(\seqs{m}\in\Kggq{m}\), then \(\mpeqKc{\seqs{m}}\cap\Kggeq{m}=\set{\Keqseq{s}{m}}\).
\ethm
\bproof
 Combine \rpropss{K.L2118}{K.L3220}.
\eproof

 Our next aim can be described as follows.
 Let \(m\in\N\) and let \(\seqs{m}\in\Kggq{m}\).
 Then an appropriate application of \rthm{A.T1514} leads us to the determination of all sequences \(\seqa{r}{m}\) which are contained in the equivalence class \(\mpeqKc{\seqs{m}}\).
 
\bpropl{K.P0830}
 Let \(m\in\N\) and let \(\seqs{m}\in\Kggq{m}\).
 Then \(\mpeqKc{\seqs{m}}\) coincides with the set of all sequences \(\seqa{r}{m}\) of complex \tqqa{matrices} fulfilling \(\ran{r_{m}-\Rku{m}}\cap\ran{\kpu{m-1}}=\set{\Ouu{q}{1}}\), \(r_{m}-\Rku{m}\in\Cggq\), and \(r_j=\su{j}\) for all \(j\in\mn{0}{m-1}\).
\eprop
\bproof
 Denote by \(\seqt{m}\) the \tKeqseq{\seqs{m}}.
 By virtue of \rprop{K.L2118}, we have then \(\seqt{m}\in\Kggeq{m}\) and \(\seqt{m}\mpeqK\seqs{m}\).
 According to \rdefn{K.D.Keqseq}, furthermore \(t_{m}=\Rku{m}\) and \(t_{j}=\su{j}\) for all \(j\in\mn{0}{m-1}\).
 In particular, \(\kpuo{m-1}{t}=\kpu{m-1}\) and \(\umgo{m-1}{t}=\umg{m-1}\).
 
 Consider an arbitrary \(\seqa{r}{m}\in\mpeqKc{\seqs{m}}\), \tie{}, \(\seqa{r}{m}\in\Kggq{m}\) with \(\seqa{r}{m}\mpeqK\seqs{m}\).
 In particular, \(\seqt{m}\mpeqK\seqa{r}{m}\), by virtue of \rrem{K.R.Keq}.
 From \rprop{K.L3220} we can conclude then that \(\seqt{m}\) coincides with the  \tKeqseq{\seqa{r}{m}}.
 In view of \rdefn{K.D.Keqseq}, consequently \(t_{m}=\Rkuo{m}{r}\) and \(t_{j}=r_j\) for all \(j\in\mn{0}{m-1}\).
 Hence, \(\Rku{m}=\Rkuo{m}{r}\) and \(\su{j}=r_j\) for all \(j\in\mn{0}{m-1}\) follow.
 In particular, \(\kpu{m-1}=\kpuo{m-1}{r}\) and \(\umg{m-1}=\umgo{m-1}{r}\).
 From \rrem{K.R.Q>}, we see furthermore \(\kpuo{m}{r}\in\Cggq\).
 Setting  \(A\defeq\kpuo{m}{r}\) and \(\cV\defeq\ran{\kpuo{m-1}{r}}\) we obtain, in view of \rnota{K.N.grG} and \rrem{K.R.Q=s-a}, then
\[
 r_{m}-\Rku{m}
 =r_{m}-\Rkuo{m}{r}
 =r_{m}-\umgo{m-1}{r}-\SO{A}{\cV}
 =\kpuo{m}{r}-\SO{A}{\cV}
 =A-\SO{A}{\cV}.
\]
 Because of \(\ran{\kpu{m-1}}=\cV\) and \rthmp{A.T1514}{A.T1514.a}, thus \(\ran{r_{m}-\Rku{m}}\cap\ran{\kpu{m-1}}=\set{\Ouu{q}{1}}\).
 Furthermore, \rpropp{sa1ab3}{sa1ab3.a} yields \(r_{m}-\Rku{m}\in\Cggq\).
 
 Conversely, consider an arbitrary sequence \(\seqa{r}{m}\) of complex \tqqa{matrices} fulfilling \(\ran{r_{m}-\Rku{m}}\cap\ran{\kpu{m-1}}=\set{\Ouu{q}{1}}\), \(r_{m}-\Rku{m}\in\Cggq\), and \(r_j=\su{j}\) for all \(j\in\mn{0}{m-1}\).
 In particular, then \(\kpuo{m-1}{r}=\kpu{m-1}\) and \(\umgo{m-1}{r}=\umg{m-1}\).
 Because of \(\seqs{m}\in\Kggq{m}\), we have \(\seqs{m-1}\in\Kggeq{m-1}\).
 Thus, \(\seqa{r}{m-1}\in\Kggeq{m-1}\).
 From \rlem{K.L1840} we infer that \(\Rku{m}\in\CHq\) and \(\umg{m-1}\lleq\Rku{m}\).
 Consequently, we can conclude \(r_{m}\in\CHq\) and \(\umg{m-1}\lleq\Rku{m}\lleq r_{m}\).
 Hence, \(\umgo{m-1}{r}\lleq r_{m}\), implying \(\kpuo{m}{r}\in\Cggq\), by virtue of \rrem{K.R.Q=s-a}.
 Using \rprop{K.P1430}, we obtain then \(r_m\in\exKgg{r}{m-1}\), \tie{}, \(\seqa{r}{m}\in\Kggq{m}\).
 Denote by $\Keqseq{r}{m}$ the \tKeqseq{\seqa{r}{m}} and let \(A\defeq\kpuo{m}{r}\) and \(\cV\defeq\ran{\kpuo{m-1}{r}}\).
 By \rdefn{K.D.Keqseq}, then \(\Keqseqmat{r}{m}=\umgo{m-1}{r}+\SO{A}{\cV}\) and \(\Keqseqmat{r}{j}=r_{j}\) for all \(j\in\mn{0}{m-1}\).
 Consequently, \(\Keqseqmat{r}{j}=r_{j}=\su{j}=t_j\) for all \(j\in\mn{0}{m-1}\).
 Setting \(X\defeq t_{m}-\umgo{m-1}{r}\) and  \(Y\defeq r_{m}-t_{m}\), we have, by virtue of \rrem{K.R.Q=s-a}, then \(X+Y=r_{m}-\umgo{m-1}{r}=\kpuo{m}{r}=A\) and \(X=t_{m}-\umg{m-1}=t_{m}-\umgo{m-1}{t}=\kpuo{m}{t}\) and furthermore \(Y=r_{m}-\Rku{m}\).
 In particular, \(\ran{Y}\cap\ran{\kpu{m-1}}=\set{\Ouu{q}{1}}\) and \(Y\in\Cggq\) by assumption.
 From \rrem{K.R.Ke<gg} we infer \(\seqt{m}\in\Kggq{m}\) and, consequently, \(t_m\in\exKgg{t}{m-1}\).
 In particular, \(\seqt{m-1}\in\Kggeq{m-1}\).
 \rthm{K.P1016} yields then \(\kpuo{m}{t}\in\Cggq\) and \(\ran{\kpuo{m}{t}}\subseteq\ran{\kpuo{m-1}{t}}\).
 Hence, \(X\in\Cggq\).
 Taking into account \(\kpuo{m-1}{t}=\kpu{m-1}=\kpuo{m-1}{r}\), we see furthermore \(\ran{X}\subseteq\cV\) and \(\ran{Y}\cap\cV=\set{\Ouu{q}{1}}\).
 From \rthmp{A.T1514}{A.T1514.b} we get then \(X=\SO{A}{\cV}\).
 Hence, \(\Keqseqmat{r}{m}=\umgo{m-1}{r}+\SO{A}{\cV}=t_{m}\) follows.
 Thus, the sequences $\Keqseq{r}{m}$  and \(\seqt{m}\) coincide.
 Using \rprop{K.L2118} and \rrem{K.R.Keq}, we get then \(\seqa{r}{m}\mpeqK\Keqseq{r}{m}\mpeqK\seqt{m}\mpeqK\seqs{m}\).
 Consequently, \(\seqa{r}{m}\in\mpeqKc{\seqs{m}}\).
\eproof

\bibliography{177arxiv}
\bibliographystyle{bababbrv}

\vfill\noindent
\begin{minipage}{0.5\textwidth}
 Universit\"at Leipzig\\
 Fakult\"at f\"ur Mathematik und Informatik\\
 PF~10~09~20\\
 D-04009~Leipzig
\end{minipage}
\begin{minipage}{0.49\textwidth}
 \begin{flushright}
  \texttt{
   fritzsche@math.uni-leipzig.de\\
   kirstein@math.uni-leipzig.de\\
   maedler@math.uni-leipzig.de
  } 
 \end{flushright}
\end{minipage}

\end{document}